\documentclass[a4paper]{article}
\usepackage[french]{RR}
\usepackage{amsmath,amssymb}

\usepackage{algorithm,algpseudocode,algorithmicx}  % For pseudocode algorithms
\usepackage{listings,url,verbatim}    % Could we do without listings?
 
\usepackage{subfigure}
\usepackage{longtable}
\usepackage{lscape}
\usepackage{multirow}
\usepackage{graphicx}
\usepackage{rotating}

\newcommand{\lt}{\left}
\newcommand{\rt}{\right}

\newcommand\abs[1]{\lvert#1\rvert}
\newcommand\norm[1]{\lVert#1\rVert}
\newcommand{\prepA}{{\bar{A}}}

\newcommand{\secA}{{\hat{A}}}

\newcommand{\algA}{{LU\_PRRP}}
\newcommand{\algB}{{CALU\_PRRP}}

\RRtitle{LU factorization with panel rank revealing pivoting and its communication avoiding version}
\RRetitle{}

\RRauthor{Amal Khabou \thanks{ Laboratoire de Recherche en Informatique, Universit\'e Paris-Sud 11, INRIA Saclay - Ile de France} ({\tt amal.khabou@inria.fr}) \and James W. Demmel\thanks{Computer Science Division and Mathematics Department, UC Berkeley, CA 94720-1776, USA. This work has been supported in part by Microsoft (Award \#024263) and Intel (Award \#024894) funding and by matching funding by U.C. Discovery (Award \#DIG07-10227). Additional support comes from Par Lab affiliates National Instruments, Nokia, NVIDIA, Oracle, and Samsung. This work has also been supported by DOE grants
DE-SC0003959, DE-SC0004938, and DE-AC02-05CH11231} ({\tt demmel@cs.berkeley.edu}) \and Laura Grigori\thanks{INRIA Saclay - Ile de France, Laboratoire de Recherche en Informatique, Universit\'e Paris-Sud 11, France. This work has been supported in part by the French National Research Agency (ANR) through COSINUS program (project PETALH no ANR-10-COSI-013).} ({\tt laura.grigori@inria.fr})  \and Ming Gu \thanks{Mathematics Department, UC Berkeley, CA 94720-1776, USA} ({\tt mgu@math.berkeley.edu}). }
\titlehead{LU\_PRRP and CALU\_PRRP}
\RRresume{On pr\' esente dans ce rapport une factorisation LU qui est plus stable que l'\' elimination de Gauss avec pivotage partiel (GEPP) en terme de
facteur de croissace. Cette factorisation permet de r\' esoudre des matrices pathologiques o\`u GEPP \' echoue (Wilkinson, Foster). On pr\' esente aussi la version
minimisant la communication de ce nouvel algorithme. }
\RRabstract{
We present the LU decomposition with panel rank revealing pivoting
(LU\_PRRP), an LU factorization algorithm based on strong rank
revealing QR panel factorization.
LU\_PRRP is more stable than Gaussian elimination with partial
pivoting (GEPP), with a theoretical upper bound of the growth factor
of $(1+ \tau b)^{n\over b}$, where $b$ is the size of the panel used
during the block factorization, $\tau$ is a parameter of the strong
rank revealing QR factorization, and $n$ is the number of columns of
the matrix. For example, if the size of the panel is $b = 64$, and
$\tau = 2$, then $(1+2b)^{n/b} = (1.079)^n \ll 2^{n-1}$, where
$2^{n-1}$ is the upper bound of the growth factor of GEPP.  Our
extensive numerical experiments show that the new factorization scheme
is as numerically stable as GEPP in practice, but it is more resistant
to pathological cases and easily solves the Wilkinson matrix and the
Foster matrix. The LU\_PRRP factorization does only $O(n^2 b)$
additional floating point operations compared to GEPP.

We also present CALU\_PRRP, a communication avoiding version of
LU\_PRRP that minimizes communication. CALU\_PRRP is based on
tournament pivoting, with the selection of the pivots at each step of
the tournament being performed via strong rank revealing QR
factorization.
CALU\_PRRP is more stable than CALU, the communication avoiding
version of GEPP, with a theoretical upper bound of the growth factor
of $(1+ \tau b)^{{n\over b}(H+1)-1}$, where $b$ is the size of the
panel used during the factorization, $\tau$ is a parameter of the
strong rank revealing QR factorization, $n$ is the number of columns of
the matrix, and $H$ is the height of the reduction tree used during
tournament pivoting.  The upper bound of the growth factor of CALU is
$2^{n(H+1)-1}$. CALU\_PRRP is also more stable in practice and is
resistant to pathological cases on which GEPP and CALU fail.
}
\RRmotcle{Factorisation LU, stabilit\'e num\'erique, minimisation de la communication, factorisation QR}
\RRkeyword{
LU factorization, numerical stability, communication avoiding, strong
rank revealing QR factorization
}

\RRprojets{Grand-Large}
\RCSaclay
\begin{document} 
\RRNo{7867}
\makeRR

\section{Introduction}
\label{sec:introduction}

The LU factorization is an important operation in numerical linear
algebra since it is widely used for solving linear systems of
equations, computing the determinant of a matrix, or as a building
block of other operations.  It consists of the decomposition of a
matrix $A$ into the product $A = \Pi LU$, where $L$ is a lower triangular
matrix, $U$ is an upper triangular matrix, and $\Pi$ a permutation
matrix.  The performance of the
LU decomposition is critical for many applications, and it has
received a significant attention over the years.  Recently large
efforts have been invested in optimizing this linear algebra kernel,
in terms of both numerical stability and performance on emerging
parallel architectures.

The LU decomposition can be computed using Gaussian elimination
with partial pivoting, a very stable operation in practice, except for
several pathological cases, such as the Wilkinson matrix
\cite{wilk61:_error_analysis, higham89:_large_growt_factor_in_gauss},
the Foster matrix \cite{foster94:_gepp_fail}, or the Wright matrix
\cite{wright92:_matrices_gepp_unstable}.  Many papers
\cite{trefethen90:_averag_gauss, skell80:_itref_gelimination,
  stewart95:_statisticge} discuss the stability of the Gaussian
elimination, and it is known
\cite{higham89:_large_growt_factor_in_gauss, gould91:_gecp,
  foster97:_ge_rookp} that the pivoting strategy used, such as
complete pivoting, partial pivoting, or rook pivoting, has an
important impact on the numerical stability of this method, which
depends on a quantity referred to as the growth factor.  However, in
terms of performance, these pivoting strategies represent a
limitation, since they require asympotically more communication than
established lower bounds on communication indicate is necessary
\cite{demmel07:_tall_skinn_qr, balard11:_min_com_linear_alg}.

Technological trends show that computing floating point operations is
becoming exponentially faster than moving data from the memory where
they are stored to the place where the computation occurs.  Due to
this, the communication becomes in many cases a dominant factor of the
runtime of an algorithm, that leads to a loss of its efficiency.  This
is a problem for both a sequential algorithm, where data needs to be
moved between different levels of the memory hierarchy, and a parallel
algorithm, where data needs to be communicated between processors.

This challenging problem has prompted research on algorithms that
reduce the communication to a minimum, while being numerically as
stable as classic algorithms, and without increasing significantly the
number of floating point operations performed
\cite{demmel07:_tall_skinn_qr, grigori2008calu}.  We refer to these
algorithms as communication avoiding.  One of the first such
algorithms is the communication avoiding LU factorization (CALU)
\cite{grigori2008calu, Grigori:EECS-2010-29}.  This algorithm is
optimal in terms of communication, that is it performs only
polylogarithmic factors more than the theoretical lower bounds on
communication require \cite{demmel07:_tall_skinn_qr,
  balard11:_min_com_linear_alg}.  Thus, it brings considerable
improvements to the performance of the LU factorization compared to
the classic routines that perform the LU decomposition such as the
PDGETRF routine of ScaLAPACK, thanks to a novel pivoting strategy
referred to as tournament pivoting. It was shown that CALU is faster
in practice than the corresponding routine PDGETRF implemented in
libraries as ScaLAPACK or vendor libraries, on both distributed
\cite{grigori2008calu} and shared memory computers
\cite{donfack2010:_multithreaded}.  While in practice CALU is as
stable as GEPP, in theory the upper bound of its growth factor is
worse than that obtained with GEPP.  One of our goals is to design an
algorithm that minimizes communication and that has a smaller upper
bound of its growth factor than CALU. 

In the first part of this paper we present the {\algA} factorization,
a novel LU decomposition algorithm based on that we call panel rank
revealing pivoting (PRRP).  The {\algA} factorization is based on a
block algorithm that computes the LU decomposition as follows.  At
each step of the block factorization, a block of columns (panel) is
factored by computing the strong rank revealing QR (RRQR)
factorization \cite{gu96:_strong_rank_revealing_qr} of its
transpose. The permutation returned by the panel rank revealing
factorization is applied on the rows of the input matrix, and the $L$
factor of the panel is computed based on the $R$ factor of the strong
RRQR factorization.  Then the trailing matrix is updated.
%The main differences with other LU factorizations as GEPP lie in the
%pivoting scheme and the computation of the $L$ factor.  However, the
%$U$ factor is computed as in other classic LU factorizations.
In exact arithmetic, the {\algA} factorization computes a block LU
decomposition based on a different pivoting strategy, the panel rank
revealing pivoting.  The factors obtained from this decomposition can
be stored in place, and so the {\algA} factorization has the same
memory requirements as standard LU and can easily replace it in any
application. 

We show that {\algA} is more stable than GEPP.  Its growth factor is
upper bounded by $(1 + \tau b)^{n \over b}$, where $b$ is the size of
the panel, $n$ is the number of columns of the input matrix, and
$\tau$ is a parameter of the panel strong RRQR factorization.  This
bound is smaller than $2^{n-1}$, the upper bound of the growth factor
for GEPP.  For example, if the size of the panel is $b = 64$, then
$(1+2b)^{n/b} = (1.079)^n \ll 2^{n-1}$.  In terms of cost, it performs
only $O(n^2 b)$ more floating point operations than GEPP.  In
addition, our extensive numerical experiments on random matrices and
on a set of special matrices show that the {\algA} factorization is
very stable in practice and leads to modest growth factors, smaller
than those obtained with GEPP.  It also solves easily pathological
cases, as the Wilkinson matrix and the Foster matrix, on which GEPP
fails.  While the Wilkinson matrix is a matrix constructed such that
GEPP has an exponential growth factor, the Foster matrix
\cite{foster97:_ge_rookp} arises from a real application.

We also discuss the backward stability of {\algA} using three metrics,
the relative error $\| PA-LU \| / \|A \|$, the normwise backward
error (\ref{eq:eta}), and the componentwise backward error (\ref{eq:comp}). For the matrices in our
set, the relative error is at most $5.26\times 10^{-14}$, the normwise
backward error is at most $1.09\times 10^{-14}$, and the componentwise
backward error is at most $3.3 \times 10^{-14}$ (with the exception of
three matrices, \textit{sprandn}, \textit{compan}, and
\textit{Demmel}, for which the componentwise backward error is $1.3
\times 10^{-13}$, $6.9 \times 10^{-12}$, and $1.16 \times 10^{-8}$
respectively).  Later in this paper, figure \ref{fig:summary} displays
the ratios of these errors versus the errors of GEPP, obtained by
dividing the maximum of the backward errors of {\algA} and the machine
epsilon $(2^{-53})$ by the maximum of those of GEPP and the machine
epsilon. For all the matrices in our set, the growth factor of {\algA}
is always smaller than that of GEPP (with the exception of one matrix,
the \textit{compar} matrix).  For random matrices, the relative error
of the factorization of {\algA} is always smaller than that of
GEPP. However, for the normwise and the componentwise backward errors,
GEPP is slightly better, with a ratio of at most $2$ between the two.
%In terms of the relative error of the factorization, the ratio is
%less than $1$ in $100\%$ of cases.
For the set of special matrices, the ratio of the relative error is at
most $1$ in over $75\%$ of cases, that is {\algA} is more stable than
GEPP. For the rest of the $25\%$ of the cases, the ratio is at most
$3$, except for one matrix (\textit{hadamard}) for which the ratio is
$23$ and the backward error is on the order of $10^{-15}$.  The ratio
of the normwise backward errors is at most $1$ in over $75\%$ of
cases, and always $3.4$ or smaller.  The ratio of the componentwise
backward errors is at most $2$ in over $81\%$ of cases, and always $3$
or smaller (except for one matrix, the \textit{compan} matrix, for
which the componentwise backward error is $6.9\times 10^{-12}$ for
{\algA} and $6.2 \times 10^{-13}$ for GEPP).

In the second part of the paper we introduce the {\algB}
factorization, the communication avoiding version of {\algA}.  It is
based on tournament pivoting, a strategy introduced in
\cite{Grigori:EECS-2010-29} in the context of CALU, a communication
avoiding version of GEPP.  With tournament pivoting, the panel
factorization is performed in two steps.  The first step selects $b$
pivot rows from the entire panel at a minimum communication cost.  For
this, sets of $b$ candidate rows are selected from blocks of the
panel, which are then combined together through a reduction-like
procedure, until a set of $b$ pivot rows are chosen.  {\algB} uses the
strong RRQR factorization to select $b$ rows at each step of the
reduction operation, while CALU is based on GEPP.  In the second step
of the panel factorization, the pivot rows are permuted to the
diagonal positions, and the QR factorization with no pivoting of the
transpose of the panel is computed.  Then the algorithm proceeds as
the {\algA} factorization.  Note that the usage of the strong RRQR
factorization ensures that bounds are respected locally at each step
of the reduction operation, but it does not ensure that the growth
factor is bounded globally as in {\algA}.

To address the numerical stability of the communication avoiding
factorization, we show that performing the {\algB} factorization of a
matrix $A$ is equivalent to performing the {\algA} factorization of a
larger matrix, formed by blocks of $A$ and zeros.  This equivalence
suggests that {\algB} will behave as {\algA} in practice and it will
be stable.  The dimension and the sparsity structure of the larger
matrix also allows us to upper bound the growth factor of {\algB} by
$(1+ \tau b)^{{n\over b}(H+1)-1}$, where in addition to the parameters
$n, b$, and $\tau$ previously defined, $H$ is the height of the
reduction tree used during tournament pivoting.

This algorithm has two significant advantages over other classic
factorization algorithms.  First, it minimizes communication, and
hence it will be more efficient than {\algA} and GEPP on architectures
where communication is expensive.  Here communication refers to both
latency and bandwidth costs of moving data between levels of the
memory hierarchy in the sequential case, and the cost of moving data
between processors in the parallel case.  Second, it is more stable
than CALU.  Theoretically, the upper bound of the growth factor of
{\algB} is smaller than that of CALU, for a reduction tree with a same
height.  More importantly, there are cases of interest for which it is
smaller than that of GEPP as well.  Given a reduction tree of height
$H=\log{P}$, where $P$ is the number of processors on which the
algorithm is executed, the panel size $b$ and the parameter $\tau$ can
be chosen such that the upper bound of the growth factor is smaller
than $2^{n-1}$.  Extensive experimental results show that {\algB} is
as stable as {\algA}, GEPP, and CALU on random matrices and a set of
special matrices.  Its growth factor is slightly smaller than that of
CALU.  In addition, it is also stable for matrices on which GEPP
fails.

As for the {\algA} factorization, we discuss the stability of {\algB}
using three metrics.  For the matrices in our set, the relative error
is at most $9.14 \times 10^{-14}$, the normwise backward error is at
most $1.37\times 10^{-14}$, and the componentwise backward error is at
most $1.14 \times 10^{-8}$ for Demmel matrix.  Figure
\ref{fig:casummary} displays the ratios of the errors with respect to
those of GEPP, obtained by dividing the maximum of the backward
errors of {\algB} and the machine epsilon by the maximum of those
of GEPP and the machine epsilon. For random matrices, all the backward
error ratios are at most $2.4$.  For the set of special matrices, the
ratios of the relative error are at most $1$ in over $62\%$ of cases,
and always smaller than $2$, except for $8\%$ of cases, where the
ratios are between $2.4$ and $24.2$. The ratios of the normwise
backward errors are at most $1$ in over $75\%$ of cases, and always
$3.9$ or smaller.  The ratios of componentwise backward errors are at
most $1$ in over $47\%$ of cases, and always $3$ or smaller, except
for $7$ ratios which have values up to $74$.

We also discuss a different version of {\algA} that minimizes
communication, but can be less stable than {\algB}, our method of
choice for reducing communication.  In this different version, the
panel factorization is performed only once, during which its
off-diagonal blocks are annihilated using a reduce-like operation,
with the strong RRQR factorization being the operator used at each
step of the reduction.  Every such factorization of a block of rows of
the panel leads to the update of a block of rows of the trailing
matrix.  Independently of the shape of the reduction tree, the upper
bound of the growth factor of this method is the same as that of
{\algA}.  This is because at every step of the algorithm, a row of the
current trailing matrix is updated only once.  We refer to the version
based on a binary reduction tree as block parallel {\algA}, and to the
version based on a flat tree as block pairwise {\algA}.  There are
similarities between these two algorithms, the LU factorization based
on block parallel pivoting (an unstable factorization), and the LU
factorization based on block pairwise pivoting (whose stability is
still under investigation)
\cite{sorensen85:_analy_of_pairw_pivot_in_gauss_elimin,
  trefethen90:_averag_gauss, barron60:_solut_of_simul_linear_equat}.
%known to be potentially unstable.
All these methods perform the panel factorization as a reduction
operation, and the factorization performed at every step of the
reduction leads to an update of the trailing matrix.  However, in
block parallel pivoting and block pairwise pivoting, GEPP is used at
every step of the reduction, and hence $U$ factors are combined
together during the reduction phase.  While in the block parallel and
block pairwise {\algA}, the reduction operates always on original rows
of the current panel.

Despite having better bounds, the block parallel {\algA} based on a
binary reduction tree of height $H = \log{P}$ is unstable for certain
values of the panel size $b$ and the number of processors $P$.  The
block pairwise {\algA} based on a flat tree of height $H = {n \over b}$
appears to be more stable.  The growth factor is larger than that of
{\algB}, but it is smaller than $n$ for the sizes of the matrices in
our test set.  Hence, potentially this version can be more stable than
block pairwise pivoting, but requires further investigation.

The remainder of the paper is organized as follows.  Section
\ref{sec:rrqrlu_algebra} presents the algebra of the {\algA}
factorization, discusses its stability, and compares it with that of
GEPP. It also presents experimental results showing that {\algA} is
more stable than GEPP in terms of worst case growth factor, and it is
more resistant to pathological matrices on which GEPP fails.  Section
\ref{sec:carrqrlu_algebra} presents the algebra of {\algB}, a
communication avoiding version of {\algA}. It describes similarities
between {\algB} and {\algA} and it discusses its stability.  The
communication optimality of {\algB} is shown in section
\ref{sec:bounds}, where we also compare its performance model with
that of the CALU algorithm.  Section \ref{sec:related_work} discusses
two alternative algorithms that can also reduce communication, but can
be less stable in practice.  Section \ref{sec:conclusion} concludes
and presents our future work.

\section{ {\algA} Method}
\label{sec:rrqrlu_algebra}

In this section we introduce the {\algA} factorization, an LU
decomposition algorithm based on panel rank revealing pivoting
strategy.  It is based on a block algorithm, that factors at each step
a block of columns (a panel), and then it updates the trailing matrix.
The main difference between {\algA} and GEPP resides in the panel
factorization.  In GEPP the panel factorization is computed using LU
with partial pivoting, while in {\algA} it is computed by performing a
strong RRQR factorization of its transpose.  This leads to a different
selection of pivot rows, and the obtained $R$ factor is used to
compute the block $L$ factor of the panel.  In exact arithmetic, {\algA}
performs a block LU decomposition with a different pivoting scheme,
which aims at improving the numerical stability of the factorization
by bounding more efficiently the growth of the elements.  We also
discuss the numerical stability of {\algA}, and we show that both in
theory and in practice, {\algA} is more stable than GEPP.

\subsection{The algebra}
{\algA} is based on a block algorithm that factors the input matrix $A$
of size $m \times n$ by traversing blocks of columns of size
$b$. Consider the first step of the factorization, with the matrix A
having the following partition,
\begin{eqnarray}
\label{eq:blockPart}
A = 
 \left[
\begin{array}{ccccc} 	
A_{11} & A_{12}\\
A_{21} & A_{22} \\
\end{array} 	
\right], 
\end{eqnarray}
where $A_{11}$ is of size $b\times b$, $A_{21}$ is of size
$(m-b)\times b$, $A_{12}$ is of size $b\times (n-b)$, and $A_{22}$ is
of size $(m-b) \times (n-b)$.

The main idea of the {\algA} factorization is to eliminate the
elements below the $b \times b$ diagonal block such that the
multipliers used during the update of the trailing matrix are bounded
by a given threshold $\tau$.  For this, we perform a strong RRQR
factorization on the transpose of the first panel of size $m\times b$
to identify a permutation matrix $\Pi$, that is $b$ pivot rows,
 $$
 \left[
\begin{array}{cc} 	
A_{11} \\
A_{21} \\
\end{array} 	
\right] ^{T} \Pi = 
\left[
\begin{array}{cc} 	
\secA_{11} \\
\secA_{21} \\
\end{array} 	
\right]^{T}
=
Q \left[  \begin{array}{cc}  R (1:b,1:b)&R(1:b,b+1:m) \end{array} 	\right] 
= Q  \left[  \begin{array}{cc} R_{11}&R_{12} \end{array} \right],
$$ where $\secA$ denotes the permuted matrix $A$.  The strong RRQR
factorization ensures that the quantity $R_{12}^T (R_{11}^{-1})^T$
is bounded by a given threshold $\tau$ in the $\max$ norm.  The strong
RRQR factorization, as described in Algorithm~\ref{srrqr} in Appendix
A, computes first the QR factorization with column pivoting, followed
by additional swaps of the columns of the $R$ factor and updates of
the QR factorization, so that $\norm{R_{12}^T (R_{11}^{-1})^T}_{\max}
\leq \tau$.  

%As explained in section \ref{sec:background}, if we note $\secA$ the
%permuted matrix, the Strong RRQR ensures that
%$\norm{\secA_{21}\secA_{11}^{-1}}_{\max} =
%\norm{R_{11}^{-1}R_{12}}_{\max} \leq \tau$.

After the panel factorization, the transpose of the computed
permutation $\Pi$ is applied on the input matrix $A$, and then the update
of the trailing matrix is performed,
\begin{eqnarray}
\label{PRRP_panel}
\secA =\Pi^T  A =
\left[
\begin{array}{cc} 
I_{b}& \\
L_{21} &I_{m - b}  \\
\end{array} 
\right] \left[
\begin{array}{cc} 
\secA_{11}& \secA_{12}\\
& {\secA}_{22}^s\\
\end{array} 
\right], 
\end{eqnarray}
where 
\begin{eqnarray}\label{update}
{\secA}_{22}^s = \secA_{22} - L_{21} \secA_{12}.
\end{eqnarray}

Note that in exact arithmetic, we have $L_{21}=\secA_{21}\secA_{11}^{-1} =
R_{12}^T (R_{11}^{-1})^T$.  Hence the factorization in equation
\eqref{PRRP_panel} is equivalent to the factorization
$$ 
\secA = \Pi^T  A =
\left[
\begin{array}{cc} 
I_{b}& \\
\secA_{21} \secA_{11}^{-1}&I_{m - b}  \\
\end{array} 
\right] \left[
\begin{array}{cc} 
\secA_{11}& \secA_{12}\\
& \secA_{22}^s\\
\end{array} 
\right], 
$$
where \begin{eqnarray}\label{update_noperm}
\secA_{22}^s= \secA_{22} - \secA_{21} \secA_{11}^{-1} \secA_{12},
\end{eqnarray}
and $\secA_{21} \secA_{11}^{-1}$ was computed in a numerically stable
way such that it is bounded in max norm by $\tau$. 

Since the block size is in general $b \geq 2$, performing {\algA} on a
given matrix $A$ first leads to a block LU factorization, with the diagonal
blocks $ \secA_{ii}$ being square of size $b \times b$.  
An additional Gaussian elimination with partial pivoting is performed
on the $b\times b$ diagonal block $\secA_{11} $ as well as the update
of the corresponding trailing matrix $\secA_{12} $. Then the
decomposition obtained after the elimination of the first panel of the
input matrix $A$ is {\small
$$
 \secA = \Pi^T A = \left[
\begin{array}{cc} 	
 I_{b}   &       \\
L_{21}  &  I_{m-b}    \\

\end{array} 	
\right] 
 \left[
\begin{array}{cc} 	
\secA_{11}     &  \secA_{12}  \\
             &   \secA_{22}^s  \\
\end{array} 
\right] 
= 
\left[
\begin{array}{cc} 	
 I_b   &       \\
L_{21}  &  I_{m-b}    \\

\end{array} 	
\right] 
\left[
\begin{array}{cc} 	
 L_{11}   &       \\
   &  I_{m-b}    \\

\end{array} 	
\right] 
 \left[
\begin{array}{cc} 	
U_{11}     &  U_{12}  \\
             &   \secA_{22}^s  \\
\end{array} 
\right],
$$ } where $L_{11}$ is a lower triangular $b\times b$ matrix with unit
diagonal and $U_{11}$ is an upper triangular $b\times b$ matrix.  We
show in section \ref{subsec:luprrpnumstab} that this step does not
affect the stability of the {\algA} factorization.  Note that the
factors $L$ and $U$ can be stored in place, and so {\algA} has the
same memory requirements as the standard LU decomposition and can
easily replace it in any application.

Algorithm \ref{rrqrlualg} presents the {\algA} factorization of a matrix
$A$ of size $n \times n$ partitioned into $n \over b$ panels.  The
number of floating-point operations performed by this algorithm is
 \begin{eqnarray*}
\# flops = { 2\over 3} n^3 + O(n^2 b),
\end{eqnarray*}
which is only $O(n^2 b)$ more floating point operations than GEPP.
The detailed counts are presented in Appendix C.  When the QR
factorization with column pivoting is sufficient to obtain the desired
bound for each panel factorization, and no additional swaps are
performed, the total cost is
 \begin{eqnarray*}
\#flops = { 2\over 3}  n^3 + { 3\over 2}n^2b.
\end{eqnarray*}

%At each step, it factors the current panel using the Strong RRQR
%factorization, and then it updates the trailing matrix as detailed in
%Algorithm \ref{rrqrlualg} where we consider the panel $ A_j$ of size
%$m-(j-1)b \times b$.
\begin{algorithm}[h!]
\caption{{\algA} factorization of a matrix A of size $n\times n$ }
\label{rrqrlualg}
\begin{algorithmic}[1]
\For{$j$ from 1 to $n\over b$}  
 \State{Let $A_j$ be the current panel $A_j = A((j-1)b+1:n,(j-1)b+1:jb)$.}
\State{Compute panel factorization $ A_j^{T} \Pi_j := Q_jR_j $ using strong RRQR factorization,}
\State{$\; \; \; L_{2j} := ( {R_j(1:b,1:b)}^{-1} R_j(1:b,b+1:n-(j-1)b))^T$.}
\State{Pivot by applying the permutation matrix $\Pi_j^T$ on the entire
  matrix, $A = \Pi_j^T A$. }
\State{Update the trailing matrix,}
\State{$\; \; \; A(jb+1:n,jb+1:n) -= L_{2j} A((j-1)b+1:jb,jb+1:n)$.  } 
 \State{Let $A_{jj}$ be the current $b\times b$ diagonal block,} 
\State{$\; \; \;A_{jj} = A((j-1)b+1:jb,(j-1)b+1:jb)$.}
\State{Compute  $A_{jj} = \Pi_{jj}L_{jj} U_{jj}$ using GEPP.}
\State{Compute $U((j-1)b+1:jb,jb+1:n) = L_{jj}^{-1}\Pi_{jj}^T A((j-1)b+1:jb,jb+1:n)$.}
  \EndFor  
  \end{algorithmic}
\end{algorithm}

\subsection{Numerical stability}
\label{subsec:luprrpnumstab}

In this section we discuss the numerical stability of the {\algA}
factorization.  The stability of an LU decomposition depends on the
growth factor.  In his backward error analysis
\cite{wilk61:_error_analysis}, Wilkinson proved that the computed
solution $\hat{x}$ of the linear system $Ax = b$, where A is of size
$n\times n$, obtained by Gaussian elimination with partial pivoting or
complete pivoting satisfies
\begin{eqnarray*}
(A+ \Delta A) \hat{x} = b, \ \ \ \ \ \norm{\Delta A} _{\infty} \leq p(n) g_W u \norm{A} _{\infty}. 
\end{eqnarray*} 
In the formula, $p(n)$ is a cubic polynomial, $u$ is the machine
precision, and $g_W$ is the growth factor defined by
\begin{center}
$ g_W = $$ max_{i,j,k}|a_{i,j}^{(k)}| \over  max_{i,j}|a_{i,j}| $,
\end{center}
where $a_{i,j}^{(k)}$ denotes the entry in position $(i,j)$ obtained
after $k$ steps of elimination.  Thus the growth factor measures the
growth of the elements during the elimination. The LU factorization is
backward stable if $g_W$ is of order O(1) (in practice the method is
stable if the growth factor is a slowly growing function of $n$).
Lemma 9.6 of \cite{higham02:_accur_and_stabil_of_numer_algor} (section
9.3) states a more general result, showing that the LU factorization
without pivoting of $A$ is backward stable if the growth factor is
small.
%The used pivoting strategy should therefore aim at keeping the growth
%factor $g_W$ small.
Wilkinson \cite{wilk61:_error_analysis} showed that for partial
pivoting, the growth factor $g_W \leq 2^{n-1}$, and this bound is
attainable.  He also showed that for complete pivoting, the upper
bound satisfies $g_W \leq n^{1/2}(2. 3^{1/2}...... n^{1/ (n-1)})^{1/
  2} \sim cn^{1/2}n^{1/4 \log n}$.  In practice the growth factors are
much smaller than the upper bounds.

In the following, we derive the upper bound of the growth factor for
the {\algA} factorization.  We use the same notation as in the
previous section and we assume without loss of generality that the
permutation matrix is the identity.  It is easy to see that the growth
factor obtained after the elimination of the first panel is bounded by
$(1+\tau b)$.  At the k-th step of the block factorization, the active
matrix $A^{(k)}$ is of size $(m-(k-1)b) \times (n-(k-1)b)$, and the
decomposition performed at this step can be written as
$$
A^{(k)} = 
 \left[
\begin{array}{ccccc} 	
A^{(k)}_{11} & A^{(k)}_{12}\\
A^{(k)}_{21} & A^{(k)}_{22} \\
\end{array} 	
\right] =
\left[
\begin{array}{cc} 
I_{b}& \\
L_{21}^{(k)} & I_{m - (k+1)b}  \\
\end{array} 
\right] \left[
\begin{array}{cc} 
A^{(k)}_{11}& A^{(k)}_{12}\\
& A^{(k)s}_{22}\\
\end{array} 
\right]. 
$$
The active matrix at the (k+1)-th step is $A^{(k+1)}_{22} =
A^{(k)s}_{22} = A^{(k)}_{22} - L_{21}^{(k)} A^{(k)}_{12}$.
Then $ \max_{i,j}|a_{i,j}^{(k+1)}|\leq \max_{i,j}|a_{i,j}^{(k)}|
(1+\tau b)$ with $\max_{i,j}|L_{21}^{(k)} (i,j)| \leq \tau $ and we have
\begin{eqnarray} \label{induction} 
{g_W}^{(k+1)} \leq  {g_W}^{(k)} (1+\tau b).
\end{eqnarray}
Induction on equation (\ref{induction}) leads to a growth factor of
{\algA} performed on the $n\over b$ panels of the matrix $A$ that
satisfies
\begin{eqnarray} 
  g_W \leq (1+\tau b)^{n/b}.
\end{eqnarray}
As explained in the algebra section, the {\algA} factorization leads
first to a block LU factorization, which is completed with additional
GEPP factorizations of the diagonal $b \times b$ blocks and updates of
corresponding blocks of rows of the trailing matrix.  These additional
factorizations lead to a growth factor bounded by $2^b$ on the
trailing blocks of rows.  Since we choose $b \ll n$, we conclude that the growth factor of the
entire factorization is still bounded by $(1+\tau b)^{n/b}$.

%Table \ref{bound} summarizes bounds for the growth factor $g_W$
%derived in this section for {\algA} and also recalls bounds for GEPP. It
%considers a matrix of size $m \times (b+1)$ for which one step of the
%{\algA} factorization is performed, and also the general case of a
%matrix of size $m \times n$. 

%\begin{table}[h]
%\small
%\centering
%\caption{Bounds of the growth factor $g_W$ obtained from
%  factoring a matrix of size $m \times (b+1)$ and a matrix of size $m \times n$ using
%  {\algA}  and GEPP.}
%\label{bound}
%\begin{tabular}{| c | c |  c | }
%\hline
%         & \multicolumn{2}{| c |}{ matrix of size $m \times (b+1)$}   \\ 
%         &  {\algA}& GEPP  \\ \hline
%$g_W$ upper bound    &  $1+ \tau b$ & $2^b$  \\ \hline
%         & \multicolumn{2}{| c |}{ matrix of size $m \times n$}   \\
%         & {\algA}  &GEPP  \\\hline
%$g_W$  upper bound   & $(1+ \tau b)^{n\over b}$ & $2^{n-1}$  \\
%\hline
%\end{tabular}
%\end{table}

The improvement of the upper bound of the growth factor of {\algA} with
respect to GEPP is illustrated in Table \ref{upper_bound}, where the
panel size varies from 8 to 128, and the parameter $\tau$ is equal to
$2$.  The worst case growth factor becomes arbitrarily smaller than for GEPP, for $b \geq 64$.

\begin{table}[h]
\small
\centering
\caption{Upper bounds of the growth factor $g_W$ obtained from factoring a
  matrix of size $m \times n$ using {\algA}  with different panel sizes
  and $\tau = 2$.}
 \label{upper_bound}
\begin{tabular}{| c | c | }
\hline
b&$g_W$\\
\hline
8&  $(1.425)^{n-1}$\\
16&  $(1.244)^{n-1}$\\
32&  $(1.139)^{n-1}$ \\
64&   $(1.078)^{n-1}$\\
128&  $(1.044)^{n-1}$ \\
\hline
\end{tabular}
\end{table}

Despite the complexity of our algorithm in pivot selection, we still
compute an LU factorization, only with different pivots. Consequently,
the rounding error analysis for LU factorization still applies (see,
for example,~\cite{demmel97:_applied_numer_linear_algeb}), which
indicates that element growth is the only factor controlling the
numerical stability of our algorithm.

\subsection{Experimental results}

We measure the stability of the {\algA} factorization experimentally
on a large set of test matrices by using several metrics, as the
growth factor, the normwise backward stability, and the componentwise
backward stability. The tests are performed in Matlab. In the tests, 
in most of the cases, the panel factorization is performed by using the QR with column
pivoting factorization instead of the strong RRQR factorization.  This
is because in practice $R_{12}^T (R_{11}^{-1})^T$ is already well
bounded after performing the RRQR factorization with column pivoting
($\norm{R_{12}^T (R_{11}^{-1})^T}_{\max}$ is rarely bigger than 3).  Hence no additional
swaps are needed to ensure that the elements are well bounded.
However, for the ill-conditionned special matrices 
(condition number $\geq 10^{14}$), to get small growth factors, we 
perform the panel factorization by using the strong RRQR factorization.
In fact, for these cases, QR with column pivoting does not ensure a 
small bound for $R_{12}^T (R_{11}^{-1})^T$.

We use a collection of matrices that includes random matrices, a set
of special matrices described in Table~\ref{Tab_allTestMatrix}, and
several pathological matrices on which Gaussian elimination with
partial pivoting fails because of large growth factors. The set of
special matrices includes ill-conditioned matrices as well as sparse
matrices. The pathological matrices considered are the Wilkinson
matrix and two matrices arising from practical applications, presented
by Foster \cite{foster94:_gepp_fail} and Wright
\cite{wright92:_matrices_gepp_unstable}, for which the growth factor
of GEPP grows exponentially. The Wilkinson matrix was constructed to
attain the upper bound of the growth factor of GEPP
\cite{wilk61:_error_analysis, higham89:_large_growt_factor_in_gauss},
and a general layout of such a matrix is
$$
 A = diag(\pm 1)\left[
\begin{array}{cccccc} 	
1&0&0&\cdots&0&1\\
-1& 1&0&...&0& 1 \\
-1& -1&1&\ddots&\vdots&\vdots\\
\vdots&\vdots&\ddots&\ddots&0&1\\
-1& -1&\cdots&-1&1 &1 \\
-1& -1&\cdots&-1&-1&1\\
\end{array} 	
\right]  \times \left[
\begin{array}{ccccc}
& & & &0\\
& T& & &\vdots \\
& & & &0 \\
0 & \cdots& & 0& \theta \\
\end{array}  
\right], 
$$ where T is an $(n-1) \times (n-1)$ non-singular upper triangular
 matrix and $ \theta = \max \abs{a_{ij} } $.  We also test a
 generalized Wilkinson matrix, the general form of
 such a matrix is
$$
 A = \left[
\begin{array}{cccccc} 	
1&0&0&\cdots&0&1\\
0& 1&0&...&0& 1 \\
& &1&\ddots&\vdots&\vdots\\
\vdots&&\ddots&\ddots&0&1\\
& & & &1 &1 \\
0& &\cdots & &0 &1\\
\end{array} 	
\right]  + T^T, 
$$ where T is an $n \times n$ upper triangular matrix with zero
 entries on the main diagonal. The matlab code of the matrix A is detailed in Appendix F.

The Foster matrix represents a concrete physical example that arises
from using the quadrature method to solve a certain Volterra integral
equation and it is of the form
$$
 A = \left[
\begin{array}{cccccc} 	
1&0&0&\cdots&0&-{1 \over c}\\
-{kh \over 2}& 1- {kh\over2}&0&...&0& -{1 \over c} \\
-{kh \over 2}& -kh&1-{kh \over 2}&\ddots&\vdots&\vdots\\
\vdots&\vdots&&\ddots&0&-{1\over c}\\
-{kh \over 2}& -kh&\cdots&-kh&1-{kh \over 2}&-{1 \over c} \\
-{kh \over 2}& -kh&\cdots&-kh&-kh&1-{1 \over c} - {kh \over 2}\\
\end{array} 	
\right]. 
$$

Wright \cite{wright92:_matrices_gepp_unstable} discusses two-point boundary value
problems for which standard solution techniques give rise to matrices with exponential 
growth factor when Gaussian elimination with partial pivoting is
used. This kind of problems arise for example from the multiple shooting algorithm.
A particular example of this problem is presented by the following
matrix, 
$$
 A = \left[
\begin{array}{ccccc} 	
I& & & &I\\
-e^{Mh}& I& & & 0 \\
 & -e^{Mh}&I& &\vdots\\
&&\ddots&\ddots&0\\
& & & -e^{Mh}&I\\
\end{array} 	
\right], 
$$
where $e^{Mh} = I + Mh + O(h^2)$. 

\begin{figure}[h!]
\caption{Growth factor $g_W$ of the {\algA} factorization of random matrices.}
\label{gw_rrqrlu}
\centering
\includegraphics[scale=0.4]{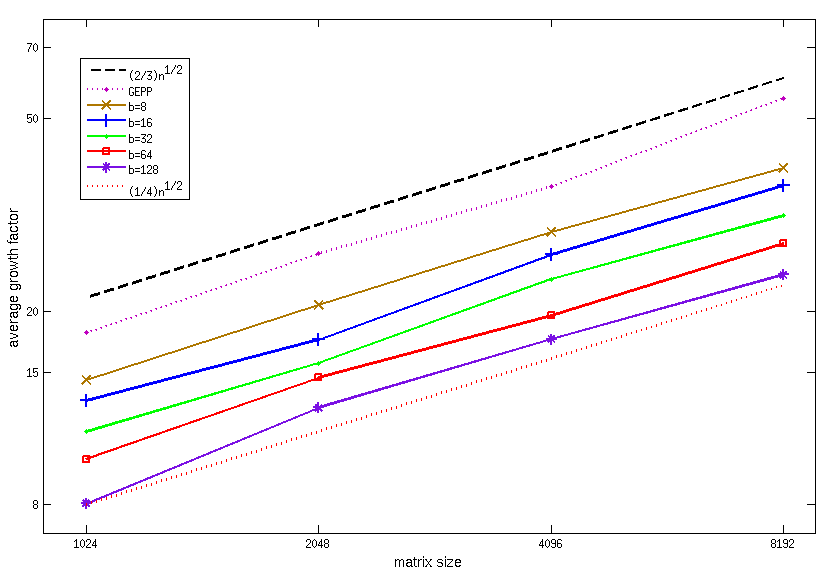}
\end{figure}

The experimental results show that the {\algA} factorization is very
stable.  Figure \ref{gw_rrqrlu} displays the growth factor of {\algA}
for random matrices of size varying from 1024 to 8192 and for sizes of
the panel varying from 8 to 128.  We observe that the smaller the size
of the panel is, the bigger the element growth is.  In fact, for a
smaller size of the panel, the number of panels and the number of
updates on the trailing matrix is bigger, and this leads to a larger
growth factor.  But for all panel sizes, the growth factor of {\algA} is
smaller than the growth factor of GEPP.  For example, for a random
matrix of size 4096 and a panel of size 64, the growth factor is only
about 19, which is smaller than the growth factor obtained by GEPP,
and as expected, much smaller than the theoretical upper bound of
$(1.078)^{4095}$.

Tables \ref{TabRRQRLU_randn1lu} and \ref{TabLUPRRP_spec1} in Appendix B
present more detailed results showing the stability of the {\algA}
factorization for random matrices and a set of special matrices.
There, we include different metrics, such as the norm of the factors,
the value of their maximum element and the backward error of the LU
factorization.
We evaluate the normwise backward stability by computing three accuracy 
tests as performed in the HPL (High-Performance Linpack) benchmark 
\cite{dongarra03:_linpac_bench}, and denoted as HPL1, HPL2 and HPL3.
\begin{eqnarray*}
\label{eq:hpl1}
\textsc{HPL1} &=&  ||Ax-b||_\infty /(\epsilon ||A||_1 *N),\\
\label{eq:hpl2}
\textsc{HPL2} &=&  ||Ax-b||_\infty /(\epsilon ||A||_1 ||x||_1),\\
\label{eq:hpl3}
\textsc{HPL3} &=&  ||Ax-b||_\infty /(\epsilon ||A||_\infty ||x||_\infty *N).
\end{eqnarray*}

In HPL, the method is considered to be accurate if the values of the
three quantities are smaller than $16$.  More generally, the values
should be of order $O(1)$. For the {\algA} factorization HPL1 is at most $8.09$,
HPL2 is at most $8.04\times 10^{-2}$ and HPL3 is at most $1.60\times 10^{-2}$. We also display the normwise backward
error, using the 1-norm,
\begin{eqnarray}
\label{eq:eta}
\eta:= \frac{||r||}{||A||~||x|| + ||b||}, 
\end{eqnarray}
and the componentwise backward error
\begin{equation}
\label{eq:comp}
w :=\max_i \frac{|r_i|}{(|A|~ | {x}| + |b|)_i},
\end{equation}
where the computed residual is $r=b-A {x}$. For our tests residuals
are computed with double-working precision.

\begin{figure}[h!]
\begin{center}
\includegraphics[angle=0,scale=0.45]{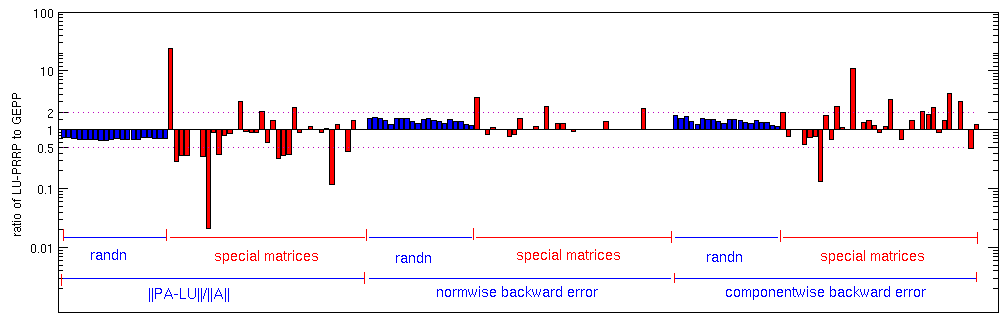}
\caption{A summary of all our experimental data, showing the ratio
  between max({\algA}'s backward error, machine epsilon) and
  max(GEPP's backward error, machine epsilon) for all the test
  matrices in our set. Each vertical bar represents such a ratio for
  one test matrix. Bars above $10^0=1$ mean that {\algA}'s backward
  error is larger, and bars below $1$ mean that GEPP's backward error
  is larger.  For each matrix and algorithm, the backward error is
  measured 3 ways. For the first third of the bars, labeled $\|PA -
  LU\| / \|A\|$, the metric is the backward error, using the Frobenius
  norm.  For the middle third of the bars, labeled ``normwise backward
  error'', the metric is $\eta$ in equation \eqref{eq:eta}. For the
  last third of the bars, labeled ``componentwise backward error'',
  the metric is $w$ in equation \eqref{eq:comp}. The test matrices are
  further labeled either as ``randn'', which are randomly generated,
  or ``special'', listed in Table \ref{Tab_allTestMatrix}.  }
\label{fig:summary}
\end{center}
\end{figure}

Figure \ref{fig:summary} summarizes all our stability results for
{\algA}.  This figure displays the ratio of the maximum between the
backward error and machine epsilon of {\algA} versus GEPP.  The
backward error is measured using three metrics, the relative error $\|
PA-LU \| / \|A \|$, the normwise backward error $\eta$, and the
componentwise backward error $w$ of {\algA} versus GEPP, and the
machine epsilon. We take the maximum of the computed error with epsilon 
since smaller values are mostly roundoff error, and so taking 
ratios can lead to extreme values with little reliability. 
Results for all the matrices in our test set are
presented, that is $20$ random matrices for which results are
presented in Table \ref{TabRRQRLU_randn1lu}, and $37$ special matrices
for which results are presented in Tables \ref{TabGEPP_spec1} and
\ref{TabLUPRRP_spec1}.  This figure shows that for random matrices,
almost all ratios are between $0.5$ and $2$.  For special matrices,
there are few outliers, up to $ 23.71$ (GEPP is more stable) for the
backward error ratio of the special matrix \textit{hadamard} and down
to $2.12\times 10^{-2}$ ({\algA} is more stable) for the backward
error ratio of the special matrix \textit{moler}.

We consider now pathological matrices on which GEPP fails.  Table
\ref{wilkinson} presents results for the linear solver using the
    {\algA} factorization for a Wilkinson matrix
    \cite{wilk85:_agebric_eigen_v_pb} of size $2048$ with a size of the
    panel varying from $8$ to $128$.  The growth factor is 1 and the
    relative error $\frac{||PA-LU||}{||A||}$ is on the order of $10 ^{-19}$.
    Table \ref{Kahan} presents results for the linear solver using the
    {\algA} algorithm for a generalized Wilkinson matrix of size $2048$
    with a size of the panel varying from $8$ to $128$.

{\small
\begin{table}[!htbp]\centering
\caption{Stability of the {\algA} factorization of a Wilkinson matrix on
  which GEPP fails.}
\label{wilkinson}
\begin{tabular}{|c|c||c|c|c|c|c|c|c|c|c|}
\hline {\sf n} & {\sf b} & $g_W$ &
$||U||_1$ & $||U^{-1}||_1$ & $||L||_1$ & $||L^{-1}||_1$ &$\frac{||PA-LU||_F}{||A||_F}$  \\\hline

\multirow{5}{*}{2048}
&    128  &  1&   1.02e+03&   6.09e+00&   1&   1.95e+00&   4.25e-20 \\ \cline{2-8}
        &    64  & 1&   1.02e+03&   6.09e+00&   1&   1.95e+00 &  5.29e-20  \\ \cline{2-8}
        &    32  & 1 &   1.02e+03&   6.09e+00&   1 &  1.95e+00&   8.63e-20 \\ \cline{2-8}
        &   16  &  1  &1.02e+03  & 6.09e+00   &1   &1.95e+00   &1.13e-19 \\ \cline{2-8}
        &    8 & 1 &   1.02e+03&   6.09e+00 &  1 & 1.95e+00 &  1.57e-19 \\ \hline
\end{tabular}
\end{table}
}

{\small
\begin{table}[!htbp]\centering
\caption{Stability of the {\algA} factorization of a generalized Wilkinson matrix on which
GEPP fails.}
\label{Kahan}
\begin{tabular}{|c|c||c|c|c|c|c|c|c|c|c|}
\hline {\sf n} & {\sf b} & $g_W$ &
$||U||_1$ & $||U^{-1}||_1$ & $||L||_1$ & $||L^{-1}||_1$ &$\frac{||PA-LU||_F}{||A||_F}$  \\\hline

\multirow{5}{*}{2048}
&    128  & 2.69&   1.23e+03&   1.39e+02&   1.21e+03&   1.17e+03 &  1.05e-15\\ \cline{2-8}
        &    64  & 2.61  & 9.09e+02   &1.12e+02  & 1.36e+03  & 1.15e+03  & 9.43e-16 \\ \cline{2-8}
        &    32  & 2.41  & 8.20e+02 &  1.28e+02&   1.39e+03 &  9.77e+02  & 5.53e-16 \\ \cline{2-8}
        &   16  & 4.08  & 1.27e+03 &  2.79e+02 &  1.41e+03  & 1.19e+03  & 7.92e-16\\ \cline{2-8}
        &    8 & 3.35 &  1.36e+03 &  2.19e+02 &  1.41e+03 &  1.73e+03 &  1.02e-15\\ \hline
\end{tabular}
\end{table}
}

For the Foster matrix, it was shown that when $c = 1$ and $kh={2 \over
  3}$, the growth factor of GEPP is $({2 \over 3})(2^{n-1} - 1)$, which
is close to the maximum theoretical growth factor of GEPP of
$2^{n-1}$.  Table \ref{tab_foster} presents results for the linear
solver using the {\algA} factorization for a Foster matrix of size
$2048$ with a size of the panel varying from $8$ to $128$ ( $c=1$,
$h=1$ and $k= {2\over 3}$ ).  According to the obtained results,
{\algA} gives a modest growth factor of $2.66$ for this practical matrix,
while GEPP has a growth factor of $10^{18}$ for the same parameters.

{\small
\begin{table}[!htbp]\centering
\caption{Stability of the {\algA} factorization of a practical matrix
  (Foster) on which GEPP fails.}
\label{tab_foster}
\begin{tabular}{|c|c||c|c|c|c|c|c|c|c|c|}
\hline {\sf n} & {\sf b} & $g_W$ &
$||U||_1$ & $||U^{-1}||_1$ & $||L||_1$ & $||L^{-1}||_1$ &$\frac{||PA-LU||_F}{||A||_F}$  \\\hline

\multirow{5}{*}{2048}
&    128  & 2.66  & 1.28e+03& 1.87e+00&1.92e+03 & 1.92e+03 &4.67e-16 \\ \cline{2-8}
        &    64  &2.66 &1.19e+03 &1.87e+00 & 1.98e+03&  1.79e+03&2.64e-16   \\ \cline{2-8}
        &    32  &2.66 &4.33e+01 & 1.87e+00&2.01e+03 &3.30e+01 &2.83e-16 \\ \cline{2-8}
        &   16  & 2.66&1.35e+03 & 1.87e+00 &  2.03e+03&2.03e+00 & 2.38e-16  \\ \cline{2-8}
        &    8 &2.66 &  1.35e+03  &  1.87e+00 &2.04e+03&2.02e+00&5.36e-17    \\ \hline
\end{tabular}
\end{table}
}

For matrices arising from the two-point boundary value problems
described by Wright, it was shown that when $h$ is chosen small enough
such that all elements of $e^{Mh}$ are less than 1 in magnitude, the growth factor obtained
using GEPP is exponential.  For
our experiment the matrix $M = \left[
\begin{array}{cc} 	
 -{1\over 6}&1\\
1&  -{1\over 6}\\
\end{array} 	
\right] $, that is $ e^{Mh} \approx \left[
\begin{array}{cc} 	
1 - {h\over 6}&h\\
h& 1 -{h\over 6}\\
\end{array} 	
\right] $, and $h=0.3$. Table \ref{tab_wright} presents results for
the linear solver using the {\algA} factorization for a Wright matrix
of size $2048$ with a size of the panel varying from $8$ to $128$.
According to the obtained results, again {\algA} gives minimum possible pivot
growth $1$ for this practical matrix, compared to the GEPP method
which leads to a growth factor of $10^{95}$ using the same parameters.
{\small
\begin{table}[!htbp]\centering
\caption{Stability of the {\algA} factorization on a practical matrix
  (Wright) on which GEPP fails.}
\label{tab_wright}
\begin{tabular}{|c|c||c|c|c|c|c|c|c|c|c|}
\hline {\sf n} & {\sf b} & $g_W$ &
$||U||_1$ & $||U^{-1}||_1$ & $||L||_1$ & $||L^{-1}||_1$ &$\frac{||PA-LU||_F}{||A||_F}$  \\\hline

\multirow{5}{*}{2048}
&    128  & 1&3.25e+00 &  8.00e+00 &  2.00e+00  & 2.00e+00 &  4.08e-17 \\ \cline{2-8}
        &    64  & 1& 3.25e+00 &  8.00e+00 &  2.00e+00 &  2.00e+00  & 4.08e-17 \\ \cline{2-8}
        &    32  & 1& 3.25e+00  & 8.00e+00  & 2.05e+00  & 2.07e+00 &  6.65e-17\\ \cline{2-8}
        &   16  & 1&  3.25e+00  & 8.00e+00 &  2.32e+00  & 2.44e+00  & 1.04e-16\\ \cline{2-8}
        &    8 & 1&  3.40e+00 &  8.00e+00 &  2.62e+00 &  3.65e+00 &  1.26e-16\\ \hline
\end{tabular}
\end{table}
}

All the previous tests show that the {\algA} factorization is very
stable for random, and for more special matrices, and it also gives modest
growth factor for the pathological matrices on which GEPP fails.  We
note that we were not able to find matrices for which {\algA} attains
the upper bound of ${\displaystyle (1+ \tau b)^{n\over b}}$ for the
growth factor.

\section{Communication avoiding {\algA}}
\label{sec:carrqrlu_algebra}

In this section we present a communication avoiding version of the
{\algA} algorithm, that is an algorithm that minimizes communication,
and so it will be more efficient than {\algA} and GEPP on architectures
where communication is expensive.  We show in this section that this
algorithm is more stable than CALU, an existing communication avoiding
algorithm for computing the LU factorization
\cite{Grigori:EECS-2010-29}.  More importantly, its parallel version
is also more stable than GEPP (under certain conditions).

\subsection{Matrix algebra}
\label{subsec:carrqrlu_algebra}
{\algB} is a block algorithm that uses tournament pivoting, a strategy
introduced in \cite{Grigori:EECS-2010-29} that allows to minimize
communication.  As in {\algA}, at each step the factorization of the
current panel is computed, and then the trailing matrix is updated.
However, in {\algB} the panel factorization is performed in two steps.
The first step, which is a preprocessing step, uses a reduction
operation to identify $b$ pivot rows with a minimum amount of
communication.  The strong RRQR factorization is the operator used at
each node of the reduction tree to select a new set of $b$ candidate
rows from the candidate rows selected at previous stages of the
reduction.  The $b$ pivot rows are permuted into the diagonal
positions, and then the QR factorization with no pivoting of the
transpose of the entire panel is computed.

In the following we illustrate tournament pivoting on the first panel,
with the input matrix $A$ partitioned as in equation
\eqref{eq:blockPart}.  Tournament pivoting considers that the first
panel is partitioned into $P=4$ blocks of rows,
\begin{eqnarray*}
A(:, 1:b)=
\left[
\begin{array}{c} 	
A_{00} \\
A_{10}\\
A_{20} \\
A_{30} \\
\end{array} 	
\right].
\end{eqnarray*}
The preprocessing step uses a binary reduction tree in our example,
and we number the levels of the reduction tree starting from $0$.  At
the leaves of the reduction tree, a set of $b$ candidate rows are
selected from each block of rows $A_{i0}$ by performing the strong
RRQR factorization on the transpose of each block $A_{i0}$.  This
gives the following decomposition, {\small
\begin{eqnarray*}
\left[
\begin{array}{c} 	
A_{00}^T \Pi_{00} \\
A_{10} ^T\Pi_{10}\\
A_{20}^T\Pi_{20} \\
A_{30}^T\Pi_{30} \\
\end{array} 	
\right]
= 
\left[
\begin{array}{c} 
Q_{00} R_{00} \\	
Q_{10} R_{10} \\
Q_{20} R_{20} \\
Q_{30} R_{30} \\
\end{array} 	
\right],   
\end{eqnarray*}
which can be written as
\begin{eqnarray*}
A(:,1:b)^T \bar{\Pi}_0= 
 A(:,1:b)^T\left[
\begin{array}{cccc} 	
 \Pi_{00} &       &       & \\ 
        & \Pi_{10} &       &  \\ 
        &       & \Pi_{20} & \\ 
        &       &       & \Pi_{30} \\
\end{array} 	
\right] 
& = &
\left[
\begin{array}{c} 
Q_{00} R_{00} \\	
Q_{10} R_{10} \\
Q_{20} R_{20} \\
Q_{30} R_{30} \\
\end{array} 	
\right],   \\
\end{eqnarray*}
}where $\bar{\Pi}_0$ is an $ m \times m$ permutation matrix with
diagonal blocks of size $ {m\over P} \times {m\over P}$, $Q_{i0}$ is
an orthogonal matrix of size $ b \times b$, and each factor $R_{i0}$
is an $ b \times {m\over P}$ upper triangular matrix.

There are now $P=4$ sets of candidate rows.  At the first level of the
binary tree, two matrices $A_{01}$ and $A_{11}$ are formed by
combining together two sets of candidate rows,
\begin{eqnarray*}
A_{01} &=&\left[\begin{array}{c} (A_{00}^T \Pi_{00})(:,1:b)\\ (A_{10}^T \Pi_{10})(:,1:b)\end{array} \right]  \\
A_{11} &=&\left[\begin{array}{c} (A_{20}^T \Pi_{20})(:,1:b)\\ (A_{30}^T \Pi_{30})(:,1:b)\end{array} \right].  \\
\end{eqnarray*}
Two new sets of candidate rows are identified by performing the strong
RRQR factorization of each matrix $A_{01}$ and $A_{11}$,
\begin{eqnarray*}
A_{01}^T \Pi_{01}= Q_{01} R_{01}, \\
A_{11}^T \Pi_{11}= Q_{11} R_{11},\\
\end{eqnarray*}
where $ \Pi_{10}$, $\Pi_{11}$ are permutation matrices of size $ 2b
\times 2b$, $Q_{01}$, $Q_{11}$ are orthogonal matrices of size $ b
\times b$, and $R_{01}, R_{11}$ are upper triangular factors of size $
b \times 2b$.

The final $b$ pivot rows are obtained by performing one last strong
RRQR factorization on the transpose of the following $ b \times 2b$
matrix :
\begin{eqnarray*}
A_{02}=\left[\begin{array}{c}(A_{01}^T \Pi_{01})(:,1:b)\\ (A_{11}^T \Pi_{11})(:,1:b)\end{array} \right],
\end{eqnarray*}
that is
\begin{eqnarray*}
A_{02}^T \Pi_{02}= Q_{02} R_{02}, \\
\end{eqnarray*}
where $\Pi_{02}$ is a permutation matrix of size $ 2b \times 2b$,
$Q_{02}$ is an orthogonal matrix of size $b \times b$, and $R_{02}$ is
an upper triangular matrix of size $ b \times 2b$.  This operation is
performed at the root of the binary reduction tree, and this ends the
first step of the panel factorization.  In the second step, the final
pivot rows identified by tournament pivoting are permuted to the
diagonal positions of $A$,
\begin{eqnarray*}
\secA =  \bar{\Pi}^T A =  \bar{\Pi}_2^T \bar{\Pi}_1^T \bar{\Pi}_0^T A,
\end{eqnarray*}
where the matrices $\bar{\Pi}_i$ are obtained by extending the
matrices $\bar{\Pi}$ to the dimension $ m\times m$, that is
\begin{eqnarray*}
\bar{\Pi}_1= \left[
\begin{array}{cc}
 	\bar{\Pi}_{01}& \\ 
& 	\bar{\Pi}_{11}\\ 
\end{array} 	
\right], 
\end{eqnarray*}
with $\bar{\Pi}_{i1}$, for $i= 0, 1$ formed as
\begin{eqnarray*}
\bar{\Pi}_{i1} = \left[
\begin{array}{cccc}
 \Pi_{i1}(1:b,1:b) &       &  \Pi_{i1}(1:b,b+1:2b)     &  \\ 
        & I_{{m\over P}-b} &       &  \\ 
  \Pi_{i1}(b+1:2b,1:b)      &       & \Pi_{i1}(b+1:2b,b+1:2b) & \\ 
     &  &  &I_{{m\over P}-b}        \\ 
\end{array} 	
\right], 
\end{eqnarray*}
and
\begin{eqnarray*}
\bar{\Pi}_2= \left[
\begin{array}{cccc} 	
 \Pi_{02}(1:b,1:b) &       &  \Pi_{02}(1:b,b+1:2b)     &  \\ 
        & I_{2{m\over P}-b} &       &  \\ 
  \Pi_{02}(b+1:2b,1:b)      &       & \Pi_{02}(b+1:2b,b+1:2b) &\\ 
        &       &       & I_{2{m\over P}-b} \\
\end{array} 	
\right]. 
\end{eqnarray*}

Once the pivot rows are in the diagonal positions, the QR
factorization with no pivoting is performed on the transpose of the
first panel,
\begin{eqnarray*}
\secA^T(1:b,:) = QR = \left[  \begin{array}{cc} R_{11}&R_{12} \end{array} \right]. 
\end{eqnarray*}
This factorization is used to update the trailing matrix, and the
elimination of the first panel leads to the following decomposition
(we use the same notation as in section \ref{sec:rrqrlu_algebra}),
$$ 
\secA  =
\left[
\begin{array}{cc} 
I_{b}& \\
\secA_{21}{\secA_{11}}^{-1}&I_{m - b}  \\
\end{array} 
\right] \left[
\begin{array}{cc} 
\secA_{11}& \secA_{12}\\
& {\secA}_{22}^{s} \\
\end{array} 
\right],
$$
where 
\begin{eqnarray*}
%\label{update_noperm_ca}
{\secA}_{22}^{s} = \secA_{22} - \secA_{21}{\secA_{11}}^{-1} \secA_{12}.
\end{eqnarray*}
As in the {\algA} factorization, the {\algB} factorization computes a block
LU factorization of the input matrix $A$. To obtain the full LU
factorization, an additional GEPP is performed on the diagonal block
$\secA_{11}$, followed by the update of the block row $\secA_{12}$.
The {\algB} factorization continues the same procedure on the trailing
matrix ${\secA}_{22}^{s}$.

Note that the factors $L$ and $U$ obtained by the {\algB} factorization
are different from the factors obtained by the {\algA} factorization.
The two algorithms use different pivot rows, and in particular the
factor $L$ of {\algB} is no longer bounded by a given threshold $\tau$
as in {\algA}. This leads to a different worst case growth factor for
{\algB}, that we will discuss in the following section.

The following figure displays the binary tree based tournament
pivoting performed on the first panel using an arrow notation (as in
\cite{Grigori:EECS-2010-29}).  The function $f(A_{ij})$ computes a
strong RRQR of the matrix $A^T_{ij}$ to select a set of $b$ candidate
rows. At each node of the reduction tree, two sets of $b$ candidate
rows are merged together and form a matrix $A_{ij}$, the function $f$ is
applied on $A_{ij}$, and another set of $b$ candidate rows is
selected.  While in this section we focused on binary trees,
tournament pivoting can use any reduction tree, and this allows the
algorithm to adapt on different architectures.  Later in the paper we
will consider also a flat reduction tree.

\begin{center}
\setlength{\unitlength}{.5cm}
\begin{picture}(7,4)

\put(0.2,0.5){$A_{30}$}
\put(0.2,1.5){$A_{20}$}
\put(0.2,2.5){$A_{10}$}
\put(0.2,3.5){$A_{00}$}

\put(1.5,0.5){$\rightarrow$}
\put(1.5,1.5){$\rightarrow$}
\put(1.5,2.5){$\rightarrow$}
\put(1.5,3.5){$\rightarrow$}

\put(2.5,0.5){$f(A_{30})$}
\put(2.5,1.5){$f(A_{20})$}
\put(2.5,2.5){$f(A_{10})$}
\put(2.5,3.5){$f(A_{00})$}

\put(4.7,0.65){$\nearrow$}
\put(4.7,1.35){$\searrow$}
\put(4.7,2.65){$\nearrow$}
\put(4.7,3.35){$\searrow$}

\put(5.5,1.0){$f(A_{11})$}
\put(5.5,3.0){$f(A_{01})$}

\put(7.7,1.5){$\nearrow$}
\put(7.7,2.5){$\searrow$}

\put(8.5,2.0){$f(A_{02})$}
\end{picture}
\end{center}

\subsection{Numerical Stability of {\algB}}

In this section we discuss the stability of the {\algB} factorization
and we identify similarities with the {\algA} factorization.  We also
discuss the growth factor of the {\algB} factorization, and we show that
its upper bound depends on the height of the reduction tree.  For the
same reduction tree, this upper bound is smaller than that obtained
with CALU.  More importantly, for cases of interest, the upper bound
of the growth factor of {\algB} is also smaller than that obtained with
GEPP.

To address the numerical stability of {\algB}, we show that performing
{\algB} on a matrix $A$ is equivalent to performing {\algA} on a larger
matrix $A_{{\algA}}$, which is formed by blocks of $A$ (sometimes
slightly perturbed) and blocks of zeros.  This reasoning is also used
in \cite{Grigori:EECS-2010-29} to show the same equivalence between
CALU and GEPP.  While this similarity is explained in detail in
\cite{Grigori:EECS-2010-29}, here we focus only on the first step of
the {\algB} factorization.  We explain the construction of the larger
matrix $A_{{\algA}}$ to expose the equivalence between the first step of
the {\algB} factorization of $A$ and the {\algA} factorization of
$A_{{\algA}}$.

Consider a nonsingular matrix $A$ of size $m \times n$ and the first
step of its {\algB} factorization using a general reduction tree of
height $H$.  Tournament pivoting selects $b$ candidate rows at each
node of the reduction tree by using the strong RRQR factorization.
Each such factorization leads to an $L$ factor which is bounded
locally by a given threshold $\tau$.  However this bound is not
guaranteed globally.  When the factorization of the first panel is
computed using the $b$ pivot rows selected by tournament pivoting, the
$L$ factor will not be bounded by $\tau$.  This results in a larger
growth factor than the one obtained with the {\algA} factorization.
Recall that in {\algA}, the strong RRQR factorization is performed on
the transpose of the whole panel, and so every entry of the obtained
lower triangular factor $L$ is bounded by $\tau$.

However, we show now that the growth factor obtained after the first
step of the {\algB} factorization is bounded by $(1+\tau b)^{H+1}$.
Consider a row $j$, and let $A^s(j, b+1:n)$ be the updated row
obtained after the first step of elimination of {\algB}.  
%To bound the growth factor, we show that $A^s(j, b+1:n)$ can be obtained by
%performing {\algA} on a larger matrix $A_{{\algA}}$, whose entries are of
%the same magnitude as the entries of the original matrix $A$.  
Suppose that row $j$ of $A$ is a candidate row at level $k-1$ of the
reduction tree, and so it participates in the strong RRQR
factorization computed at a node $s_k$ at level $k$ of the reduction
tree, but it is not selected as a candidate row by this factorization.
We refer to the matrix formed by the candidate rows at node $s_k$ as
$\prepA_k$.  Hence, row $j$ is not used to form the matrix $\prepA_k$.
Similarly, for every node $i$ on the path from node $s_k$ to the root
of the reduction tree of height $H$, we refer to the matrix formed by
the candidate rows selected by strong RRQR as $\prepA_i$.  Note that
in practice it can happen that one of the blocks of the panel is
singular, while the entire panel is nonsingular.  In this case strong
RRQR will select less than $b$ linearly independent rows that will be
passed along the reduction tree.  However, for simplicity, we assume
in the following that the matrices $\prepA_i$ are nonsingular.  For a
more general solution, the reader can consult
\cite{Grigori:EECS-2010-29}.

Let $\Pi$ be the permutation returned by the tournament pivoting
strategy performed on the first panel, that is the permutation that
puts the matrix $\prepA_H$ on diagonal.  The following equation is
satisfied,
{\small
\begin{eqnarray}
\label{eq:caluprrpequiv}
\begin{pmatrix}
\prepA_{H} &  \secA_{H} \\
A(j, 1:b) & A(j, b+1:n) \\
\end{pmatrix}
=
\begin{pmatrix}
I_{b} &  \\
A(j, 1:b)\prepA_{H}^{-1} & 1 \\
\end{pmatrix}
\cdot
\begin{pmatrix}
\prepA_{H} &  \secA_{H} \\
   & A^s(j, b+1:n) 
\end{pmatrix},
\end{eqnarray}
}
where 

\begin{eqnarray*}
\prepA_{H} &=& (\Pi A) (1:b, 1:b), \\
\secA_{H} &=& (\Pi A) (1:b, b+1:n).
\end{eqnarray*}

The updated row $A^s(j, b+1:n)$ can be also obtained by performing
{\algA} on
%the $(H-k+1)$ block columns of
a larger matrix $A_{{\algA}}$ of dimension $((H-k+1) b + 1) \times
((H-k+1) b+1) $, 
{\small
\begin{eqnarray}
\label{eq:aluprrp}
A_{{\algA}} &= &
\begin{pmatrix}
\prepA_{H} & & & & & \secA_{H} \\
\prepA_{ H-1} & \prepA_{H-1} & & & &  \\
 & \prepA_{H-2} & \prepA_{H-2} & & &  \\
 & & \ddots & \ddots & &  \\
 & & & \prepA_{k} & \prepA_{k} &  \\
 & & & & (-1)^{H-k} A(j, 1:b) & A(j, b+1:n) \\
\end{pmatrix}  \nonumber \\
& = &
\begin{pmatrix}
I_{b} & & & & &  \\
\prepA_{ H-1} \prepA_{H}^{-1} & I_{b}& & & &  \\
 & \prepA_{H-2} \prepA_{H-1}^{-1} & I_{b} & & &  \\
 & & \ddots & \ddots & &  \\
 & & & \prepA_{k} \prepA_{k+1}^{-1} &I_{b}&  \\
 & & & & (-1)^{H-k} A(j, 1:b) \prepA_{k}^{-1} & 1 \\
\end{pmatrix}
 \nonumber \\
& \cdot &
\begin{pmatrix}
\prepA_{H} & & & & & \secA_{H} \\
 & \prepA_{H-1} & & & & \secA_{H-1} \\
 &  & \prepA_{H-2} & & & \secA_{H-2} \\
 & &  & \ddots & &  \vdots \\
 & & &  & \prepA_{k} & \secA_{k} \\
 & & & &  & A^s(j, b+1:n) \\
\end{pmatrix},
\end{eqnarray}
}
where
{\small
\begin{eqnarray}
\secA_{H-i} & = & \left\{ \begin{array}{ll}
     \prepA_{H} 
             & {\rm if} \; i=0, \\
     - \prepA_{H-i} \prepA_{H-i+1}^{-1}  \secA_{H-i+1}
             & {\rm if} \; 0 < i \leq H- k. \\
             \end{array} \right.
\end{eqnarray}
}
Equation \eqref{eq:aluprrp} can be easily verified, since
\begin{eqnarray*}
A^s(j, b+1:n) &=& A(j, b+1:n) - (-1)^{H-k} A(j, 1:b)\prepA_{k}^{-1}(-1)^{H-k}\secA_{k}\\
   &=& A(j, b+1:n) - A(j, 1:b)\prepA_{k}^{-1}\prepA_{k} \prepA_{k+1}^{-1} \hdots \prepA_{H-2} \prepA_{H-1}^{-1}\prepA_{ H-1} \prepA_{H}^{-1}\secA_{H}\\
   &=& A(j, b+1:n) - A(j, 1:b) \prepA_H^{-1} \secA_{H}. \\
\end{eqnarray*}

Equations \eqref{eq:caluprrpequiv} and \eqref{eq:aluprrp} show that
the Schur complement obtained after each step of performing the
CALU\_PRRP factorization of a matrix $A$ is equivalent to the Schur
complement obtained after performing the LU\_PRRP factorization of a
larger matrix $A_{LU\_PRRP}$, formed by blocks of $A$ (sometimes
slightly perturbed) and blocks of zeros.  More generally, this implies
that the entire {\algB} factorization of $A$ is equivalent to the
{\algA} factorization of a larger and very sparse matrix, formed by
blocks of $A$ and blocks of zeros (we omit the proofs here, since they
are similar with the proofs presented in \cite{Grigori:EECS-2010-29}).

Equation \eqref{eq:aluprrp} is used to derive the upper bound of the
growth factor of {\algB} from the upper bound of the growth factor of
{\algA}.  The elimination of each row of the first panel using {\algB} can
be obtained by performing {\algA} on a matrix of maximum dimension $m
\times b(H+1)$.  Hence the upper bound of the growth factor obtained
after one step of {\algB} is $(1+\tau b)^{H+1}$.  This leads to an upper
bound of $(1 + \tau b )^{{n \over b}(H+1)-1}$ for a matrix of size
$m\times n$.

Table \ref{caluprrp_bound} summarizes the bounds of the growth factor
of {\algB} derived in this section, and also recalls the bounds of
{\algA}, GEPP, and CALU the communication avoiding version of GEPP.  It
considers the growth factor obtained after the elimination of $b$
columns of a matrix of size $m \times (b+1)$, and also the general
case of a matrix of size $m \times n$.  As discussed in section
\ref{sec:rrqrlu_algebra} already, {\algA} is more stable than GEPP in
terms of worst case growth factor.  From Table \ref{caluprrp_bound},
it can be seen that for a reduction tree of a same height, {\algB} is
more stable than CALU.

In the following we show that {\algB} can be more stable than GEPP in
terms of worst case growth factor. 
Consider a parallel version of {\algB} based on a binary reduction tree
of height $H=\log(P)$, where $P$ is the number of processors.  The
upper bound of the growth factor becomes $(1+ \tau
b)^{{{n(logP+1)}\over{b}}-1}$, which is smaller than $2^{n(logP+1) -
  1}$, the upper bound of the growth factor of CALU.  For example, if
the threshold is $ \tau = 2$, the panel size is $b = 64$, and the
number of processors is $ P = 128 =2^7$, then ${g_W}_{ {\algB}}=
(1.7)^n$.  This quantity is much smaller than $2^{7n}$ the upper bound
of CALU, and even smaller than the worst case growth factor of GEPP of
$2^{n-1}$.  In general, the upper bound of {\algB} can be smaller than
the one of GEPP, if the different parameters $\tau$, $H$, and $b$ are
chosen such that the condition
\begin{eqnarray}
\label{eq:grfactorrel}
H \leq {b \over (\log b + \log \tau)}
\end{eqnarray}
is satisfied.  For a binary tree of height $H=\log P$, it
becomes $\log P \leq {b \over (\log b + \log \tau)}$.  This is a
condition which can be satisfied in practice, by choosing $b$ and
$\tau$ appropriately for a given number of processors $P$.  For
example, when $P \leq 512$, $ b= 64$, and $\tau = 2$, the condition
\eqref{eq:grfactorrel} is satisfied, and the worst case growth factor
of {\algB} is smaller than the one of GEPP.

However, for a sequential version of {\algB} using a flat tree of height
$H = n/b$, the condition to be satisfied becomes $n \leq {b^2\over
  (\log b + \log \tau)}$, which is more restrictive.  In practice, the
size of $b$ is chosen depending on the size of the memory, and it
might be the case that it will not satisfy the condition in equation
\eqref{eq:grfactorrel}.

\begin{table}[h]
\small
\centering
\caption{Bounds for the growth factor $g_W$ obtained from factoring a
  matrix of size $m \times (b+1)$ and a matrix of size $m \times n$
  using {{\algB}}, {{\algA}}, CALU, and GEPP.  {{\algB}} and CALU use a
  reduction tree of height $H$.  The strong RRQR used in {\algA} and
  {\algB} is based on a threshold $\tau$.  For the matrix of size $m
  \times (b+1)$, the result corresponds to the growth factor obtained
  after eliminating $b$ columns.}
\begin{tabular}{| c | c | c|c | c|}
\hline
         & \multicolumn{4}{| c|}{ matrix of size $m \times (b+1)$}   \\ 
         & TSLU(b,H) & {{\algA}}(b,H)  & GEPP &    {\algA} \\\hline
$g_W$  upper bound  &$2^{b(H+1)}$ &  $(1+ \tau b)^{{H+1}}$  & $2^b$&  $1+ \tau b$  \\ \hline
         & \multicolumn{4}{| c|}{ matrix of size $m \times n$}   \\
         & CALU & {\algB}  & GEPP &    {\algA}  \\\hline
$g_W$  upper bound &$2^{n(H+1) - 1}$ & $(1+ \tau b)^{{{n(H+1)}\over{b}}-1}$ & $2^{n-1}$ & $(1+ \tau b)^{n\over b}$ \\\hline
\end{tabular}
\label{caluprrp_bound}
\end{table}

\subsection{Experimental results}
In this section we present experimental results and show that CALU\_PRRP is stable in
practice and compare them with those obtained from CALU and GEPP in \cite{Grigori:EECS-2010-29}.  We present results 
for both the binary tree scheme and the flat tree scheme.  

As in section \ref{sec:rrqrlu_algebra}, we perform our tests on matrices whose elements follow a normal
distribution.  In {\sc Matlab} notation, the test matrix is $A = {\sf
  randn}(n,n)$, and the right hand side is $ b={\sf randn}(n,1)$.  The
size of the matrix is chosen such that $n$ is a power of $2$, that is
$n=2^k$, and the sample size is 10 if $k < 13$ and 3 if $k\geq 13$. To measure the stability of CALU\_PRRP, we
discuss several metrics, that concern the LU decomposition and the
linear solver using it, such as the growth factor, normwise and
componentwise backward errors. We also perform tests on several special matrices including sparse
matrices, they are described in Appendix B.

Figure~\ref{fig:bcaluprrp} displays the values of the growth factor
$g_W$ of the binary tree based CALU\_PRRP, for different block sizes $b$ and
different number of processors $P$.  As explained in section
\ref{subsec:carrqrlu_algebra}, the block size determines the size of the
panel, while the number of processors determines the number of block
rows in which the panel is partitioned.  This corresponds to the
number of leaves of the binary tree.  We observe that the growth
factor of binary tree based CALU\_PRRP is in the most of the cases better than GEPP. 
The curves of the growth factor lie between ${1\over 2} n^{1/2}$ and ${3\over4}
n^{1/2}$ in our tests on random matrices.
These results show that binary tree based CALU\_PRRP is stable and the growth factor values obtained 
for the different layouts are better than those obtained with binary tree based CALU.
The figure~\ref{fig:bcaluprrp} includes also the growth factor of the {\algA} method with
a panel of size $ b=64$. We note that that results are better than those of binary
tree based {\algB}.
\begin{figure}[htbp]
\begin{center}
\includegraphics[angle=0,scale=0.4]{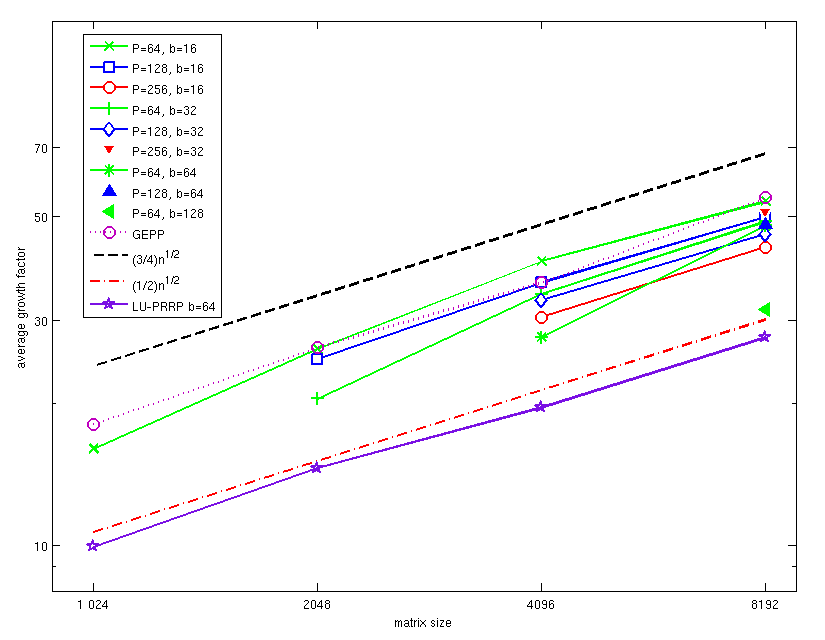}
\caption{Growth factor $g_W$ of binary tree based CALU\_PRRP for random matrices. 
} \label{fig:bcaluprrp}
\end{center}
\end{figure}

\begin{figure}[htbp]
\begin{center}
\includegraphics[angle=0,scale=0.4]{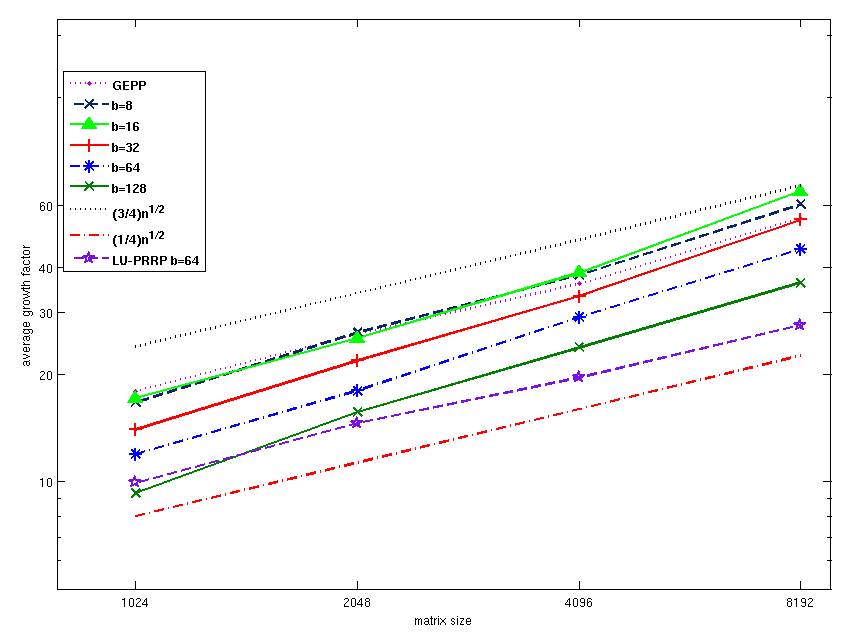}
\caption{Growth factor $g_W$ of flat tree based CALU\_PRRP for random matrices. 
} \label{fig:fcaluprrp}
\end{center}
\end{figure}

Figure~\ref{fig:fcaluprrp} displays the values of the growth factor
$g_W$ for flat tree based CALU\_PRRP with a block size $b$ varying from $8$ to
$128$. The growth factor $g_W$ is decreasing with increasing the panel size b. 
We note that the curves of the growth factor lie between  ${1\over 4} n^{1/2}$ and ${3\over 4}
n^{1/2}$ in our tests on random matrices. We also note that the results obtained with the
{\algA} method with a panel of size $b =64$ are better than those of flat tree based {\algB}.\\

The growth factors of both binary tree based and flat tree based CALU\_PRRP have 
similar (sometimes better) behavior than the growth factors of GEPP.

Table~\ref{bcaluprrp_gepp} in Appendix B presents results for the linear solver using
binary tree based \algB, together with binary tree based CALU and GEPP for
comparaison and Table~\ref{fcaluprrp_gepp} in Appendix B presents results for the linear 
solver using flat tree based \algB, together with flat tree based CALU and 
GEPP for comparaison. 
We note that for the binary tree based CALU\_PRRP, when ${m\over P} = b$, for the
 algorithm we only use $P1 = {m\over (b+1)}$ processors, since to perform a Strong 
RRQR on a given block, the number of its rows should be at least the number of its
 columns $+ 1$.  
Tables \ref{bcaluprrp_gepp} and \ref{fcaluprrp_gepp} also include results obtained 
by iterative refinement used to improve the accuracy of the solution.  For this, the
componentwise backward error in equation (\ref{eq:comp}) is used. In the previous tables, $w_b$
denotes the componentwise backward error before iterative refinement
and $N_{IR}$ denotes the number of steps of iterative refinement.
$N_{IR}$ is not always an integer since it represents an average. 
For all the matrices tested CALU\_PRRP leads to results as accurate as the results 
obtained with CALU and GEPP.

In Appendix B we present more detailed results. There we include 
some other metrics, such as the norm of the factors, the norm of
the inverse of the factors, their conditioning, the value of their 
maximum element, and the backward error of the LU factorization. 
Through the results detailed in this section and in Appendix B we 
show that binary tree based and flat tree based CALU\_PRRP are 
stable, have the same behavior as GEPP for random matrices, and 
are more stable than binary tree based and flat tree based CALU 
in terms of growth factor.

\begin{figure}[h!]
\begin{center}
\includegraphics[angle=0,scale=0.45]{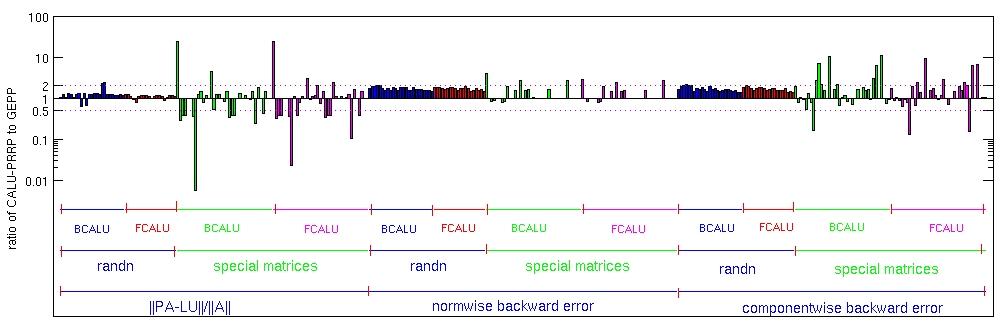}
\caption{A summary of all our experimental data, showing the ratio of
  max({\algB}'s backward error, machine epsilon) to max(GEPP's
  backward error, machine epsilon) for all the test matrices in our
  set. Each vertical bar represents such a ratio for one test
  matrix. Bars above $10^0=1$ mean that {\algB}'s backward error is
  larger, and bars below $1$ mean that GEPP's backward error is
  larger.  For each matrix and algorithm, the backward error is
  measured using three different metrics. For the last third of the
  bars, labeled ``componentwise backward error'', the metric is $w$ in
  equation \eqref{eq:comp}. The test matrices are further labeled
  either as ``randn'', which are randomly generated, or ``special'',
  listed in Table \ref{Tab_allTestMatrix}.  Finally, each test matrix
  is factored using {\algB} with a binary reduction tree (labeled
  BCALU for BCALU\_PRRP) and with a flat reduction tree (labeled FCALU
  for FCALU\_PRRP).}
\label{fig:casummary}
\end{center}
\end{figure}

Figure \ref{fig:casummary} summarizes all our stability results for
the {\algB} factorization based on both binary tree and flat tree
schemes.  As figure \ref{fig:summary}, this figure displays the ratio
of the maximum between the backward error and machine epsilon of
{\algA} versus GEPP. The backward error is measured as the relative
error $\| PA-LU \| / \|A \|$, the normwise backward error $\eta$, and
the componentwise backward error $w$.  Results for all the matrices in
our test set are presented, that is $25$ random matrices for binary
tree base {\algB} from Table \ref{Tab_bcarrqr_randn1lu}, $20$ random
matrices for flat tree based {\algB} from Table
\ref{Tab_fcarrqr_randn1lu}, and $37$ special matrices from Tables
\ref{TabFcaluprrp_spec} and \ref{TabBcaluprrp_spec}. As it can be
seen, nearly all ratios are between $0.5$ and $2.5$ for random
matrices.  However there are few outliers, for example the relative
error ratio has values between $24.2$ for the special matrix
\textit{hadamard} (GEPP is more stable than binary tree based
       {\algB}), and $5.8\times 10^{-3}$ for the special matrix
       \textit{moler} (binary tree based {\algB} is more stable than
       GEPP).
%There are also some outliers for the other ratios for some special
%matrices.

We consider now the same set of pathological matrices as in section 
\ref{sec:rrqrlu_algebra} on which GEPP fails. For the Wilkinson matrix,
 both CALU and {\algB} based on flat and binary tree give modest 
element growth. For the generalized Wilkinson matrix, the Foster matrix, 
and Wright matrix, CALU fails with both flat tree and binary tree reduction schemes. 

Tables \ref{ftpatho_lu} and \ref{btpatho_lu} present the results
obtained for the linear solver using the {\algB} factorization based
on flat and binary tree schemes for a generalized Wilkinson matrix of
size $2048$ with a size of the panel varying from $8$ to $128$ for the
flat tree scheme and a number of processors varying from $128$ to $32$
for the binary tree scheme.  The growth factor is of order 1 and the
quantity $\frac{||PA-LU||}{||A||}$ is on the order of $10
^{-16}$. Both flat tree based CALU and binary tree based CALU fail on
this pathological matrix. In fact for a generalized Wilkinson matrix
of size $1024$ and a panel of size $b = 128$, the growth factor
obtained with flat tree based CALU is of size $10^{234}$.

{\small
\begin{table}[!htbp]\centering
\caption{Stability of the flat tree based {\algB} factorization of a generalized Wilkinson matrix on which
GEPP fails.}
\label{ftpatho_lu}
\begin{tabular}{|c|c||c|c|c|c|c|c|c|c|c|}
\hline {\sf n} & {\sf b} & $g_W$ &
$||U||_1$ & $||U^{-1}||_1$ & $||L||_1$ & $||L^{-1}||_1$ &$\frac{||PA-LU||_F}{||A||_F}$  \\\hline

\multirow{5}{*}{2048}
&    128  & 2.01 &  1.01e+03 &  1.40e+02 &  1.31e+03 &  9.76e+02 &  9.56e-16\\ \cline{2-8}
        &    64  &  2.02&  1.18e+03 &  1.64e+02 &  1.27e+03 &  1.16e+03  & 1.01e-15 \\ \cline{2-8}
        &    32  & 2.04 &  8.34e+02 &  1.60e+02 &  1.30e+03 &  7.44e+02 &  7.91e-16 \\ \cline{2-8}
        &   16  &2.15  & 9.10e+02  & 1.45e+02 &  1.31e+03 &  8.22e+02 &  8.07e-16 \\ \cline{2-8}
        &    8 & 2.15  & 8.71e+02 &  1.57e+02 &  1.371e+03 &  5.46e+02 &  6.09e-16\\ \hline
\end{tabular}
\end{table}
}

{\small
\begin{table}[!htbp]\centering
\caption{Stability of the binary tree based {\algB} factorization of a generalized Wilkinson matrix on which
GEPP fails.}
\label{btpatho_lu}
\begin{tabular}{|c|c|c||c|c|c|c|c|c|c|c|c|}
\hline {\sf n} & {\sf P}& {\sf b} & $g_W$ &$||U||_1$ & $||U^{-1}||_1$ & $||L||_1$ & $||L^{-1}||_1$ &$\frac{||PA-LU||_F}{||A||_F}$  \\\hline

\multirow{5}{*}{2048}
&    128& 8  &2.10e+00 & 1.33e+03 &  1.29e+02 &  1.34e+03 &  1.33e+03 &  1.08e-15 \\ \cline{2-9}
& \multirow{2}{*}{64}   
&      16 &  2.04e+00 &  6.85e+02 &  1.30e+02 &  1.30e+03 &  6.85e+02  & 7.85e-16\\ \cline{3-9}
 &  &      8 & 8.78e+01 &  1.21e+03 &  1.60e+02  & 1.33e+03 &  1.01e+03  & 9.54e-16  \\ \cline{2-9}
& \multirow{3}{*}{32} 
        &    32  &2.08e+00  & 9.47e+02  & 1.58e+02 &  1.41e+03  & 9.36e+02 &  5.95e-16 \\ \cline{3-9}
&      &   16  & 2.08e+00 &  1.24e+03 &  1.32e+02 &  1.35e+03 &  1.24e+03  & 1.01e-15\\ \cline{3-9}
&      &    8 & 1.45e+02 &  1.03e+03 &  1.54e+02 &  1.37e+03 &  6.61e+02&   6.91e-16\\ \hline
\end{tabular}
\end{table}
}

For the Foster matrix, we have seen in the section \ref{sec:rrqrlu_algebra} that 
{\algA} gives modest pivot growth, whereas GEPP fails. 
Both flat tree based CALU and binary tree based CALU fail on the Foster matrix.
However flat tree based {\algB} and binary tree 
based {\algB} solve easily this pathological matrix. 

Tables \ref{ftfoster} and \ref{btfoster} present results for the linear
solver using the {\algB} factorization based on the flat tree scheme
and the binary tree scheme, respectively. We test a Foster matrix of size $2048$ with a panel size 
varying from $8$ to $128$ for the flat tree based {\algB} and a number of 
processors varying from $128$ to $32$ for the binary tree based {\algB}. We use the
same parameters as in section \ref{sec:rrqrlu_algebra}, that is $c=1$, $h=1$ and 
$k= {2\over 3}$. According to the obtained results, {\algB} gives a modest 
growth factor of $1.33$ for this practical matrix.

{\small
\begin{table}[!htbp]\centering
\caption{Stability of the flat tree based {\algB} factorization of a practical matrix (Foster) on which
GEPP fails.}
\label{ftfoster}
\begin{tabular}{|c|c||c|c|c|c|c|c|c|c|c|}
\hline {\sf n} & {\sf b} & $g_W$ &
$||U||_1$ & $||U^{-1}||_1$ & $||L||_1$ & $||L^{-1}||_1$ &$\frac{||PA-LU||_F}{||A||_F}$  \\\hline

\multirow{5}{*}{2048}
&    128  & 1.33&  1.71e+02  & 1.87e+00 &  1.92e+03 &  1.29e+02 &  6.51e-17\\ \cline{2-8}
        &    64  &   1.33&   8.60e+01  & 1.87e+00  & 1.98e+03  & 6.50e+01  & 4.87e-17\\ \cline{2-8}
        &    32  &     1.33&   4.33e+01 &  1.87e+00 &  2.01e+03&   3.30e+01 &  2.91e-17 \\ \cline{2-8}
        &   16  & 1.33 &  2.20e+01 &  1.87e+00  & 2.03e+03 &  1.70e+01 &  4.80e-17\\ \cline{2-8}
        &    8 & 1.33&  1.13e+01&   1.87e+00&   2.04e+03 &  9.00e+00  & 6.07e-17\\ \hline
\end{tabular}
\end{table}
}

{\small
\begin{table}[!htbp]\centering
\caption{Stability of the binary tree based {\algB} factorization of a practical matrix (Foster) on which
GEPP fails.}
\label{btfoster}
\begin{tabular}{|c|c|c||c|c|c|c|c|c|c|c|c|}
\hline {\sf n} & {\sf P}& {\sf b} & $g_W$ &$||U||_1$ & $||U^{-1}||_1$ & $||L||_1$ & $||L^{-1}||_1$ &$\frac{||PA-LU||_F}{||A||_F}$  \\\hline

\multirow{5}{*}{2048}
&    128& 8  & 1.33 &  1.13e+01 &  1.87e+00 &  2.04e+03 &  9.00e+00 &  6.07e-17\\ \cline{2-9}
& \multirow{2}{*}{64}   
&      16 & 1.33 &  2.20e+01  & 1.87e+00 &  2.03e+03 &  1.70e+01 &  4.80e-17 \\ \cline{3-9}
 &  &      8 & 1.33 & 1.13e+01  & 1.87e+00  & 2.04e+03  & 9.00e+00 &  6.07e-17 \\ \cline{2-9}
& \multirow{3}{*}{32} 
        &    32  & 1.33 & 4.33e+01 &  1.87e+00 &  2.01e+03  & 3.300e+01 &  2.91e-17 \\ \cline{3-9}
&      &   16  & 1.33 & 2.20e+01 &  1.87e+00 &  2.03e+03  & 1.70e+01  & 4.80e-17\\ \cline{3-9}
&      &    8 &  1.33&   1.13e+01 &  1.87e+00&   2.04e+03 &  9.00e+00 &  6.07e-17\\ \hline
\end{tabular}
\end{table}
}

As GEPP, both the flat tree based and the binary tree based CALU fail on the Wright matrix.
In fact for a matrix of size $2048$, a parameter $h = 0.3$, with a panel of size $b = 128$, the
flat tree based CALU gives a growth factor of $10^{98}$. With a number of processors $P = 64$ and a 
panel of size $b=16$, the binary tree based CALU also gives a growth factor of $10^{98}$. 
Tables \ref{ftwright} and \ref{btwright} present results for
the linear solver using the {\algB} factorization for a Wright matrix
of size $2048$. For the flat tree based {\algB}, the size of the panel is varying from $8$ to $128$.
For the binary tree based {\algB}, the number of processors is varying from $32$ to $128$ and the size of
the panel from $16$ to $64$ such that the number of rows in the leaf nodes is equal or bigger than two times
the size of the panel.
The obtained results, show that {\algB} gives a modest growth factor of 
$1$ for this practical matrix, compared to the CALU method.
{\small
\begin{table}[!htbp]\centering
\caption{Stability of the flat tree based {\algB} factorization of a practical matrix (Wright) on which
GEPP fails.}
\label{ftwright}
\begin{tabular}{|c|c||c|c|c|c|c|c|c|c|c|}
\hline {\sf n} & {\sf b} & $g_W$ &
$||U||_1$ & $||U^{-1}||_1$ & $||L||_1$ & $||L^{-1}||_1$ &$\frac{||PA-LU||_F}{||A||_F}$  \\\hline

\multirow{5}{*}{2048}
&    128  & 1&  3.25e+00 &  8.00e+00  & 2.00e+00 &  2.00e+00 &  4.08e-17\\ \cline{2-8}
        &    64  & 1 &  3.25e+00 &  8.00e+00 &  2.00e+00 &  2.00e+00  & 4.08e-17 \\ \cline{2-8}
        &    32  &   1 &   3.25e+00 &  8.00e+00 &  2.05e+00 &  2.02e+00 &  6.65e-17\\ \cline{2-8}
        &   16  & 1& 3.25e+00 &  8.00e+00 &  2.32e+00 &  2.18e+00 &  1.04e-16\\ \cline{2-8}
        &    8 &  1  & 3.40e+00 &  8.00e+00  & 2.62e+00  & 2.47e+00  & 1.26e-16\\ \hline
\end{tabular}
\end{table}
}
{\small
\begin{table}[!htbp]\centering
\caption{Stability of the binary tree based {\algB} factorization of a practical matrix (Wright) on which
GEPP fails.}
\label{btwright}
\begin{tabular}{|c|c|c||c|c|c|c|c|c|c|c|c|}
\hline {\sf n} & {\sf P}& {\sf b} & $g_W$ &$||U||_1$ & $||U^{-1}||_1$ & $||L||_1$ & $||L^{-1}||_1$ &$\frac{||PA-LU||_F}{||A||_F}$  \\\hline

\multirow{5}{*}{2048}
&    128& 8  & 1 & 3.40e+00 &  8.00e+00  & 2.62e+00 &  2.47e+00&   1.26e-16\\ \cline{2-9}
& \multirow{2}{*}{64}   
&      16 & 1&  3.25e+00 &  8.00e+00 &  2.32e+00  & 2.18e+00 &  1.04e-16 \\ \cline{3-9}
 &  &      8 & 1 & 3.40e+00  & 8.00e+00  & 2.62e+00 &  2.47e+00 &  1.26e-16 \\ \cline{2-9}
& \multirow{3}{*}{32} 
        &    32  &1&   3.25e+00 &  8.00e+00 &  2.05e+00 &  2.02e+00 &  6.65e-17  \\ \cline{3-9}
&      &   16  & 1  & 3.25e+00  & 8.00e+00  & 2.32e+00  & 2.18e+00  & 1.04e-16\\ \cline{3-9}
&      &    8 & 1&  3.40e+00&   8.00e+00&   2.62e+00  & 2.47e+00 &  1.26e-16\\ \hline
\end{tabular}
\end{table}
}

All the previous tests show that the {\algB} factorization is very
stable for random and more special matrices, and it also gives modest
growth factor for the pathological matrices on which CALU fails, this is 
for both binary tree and flat tree based {\algB}.  

\section{Lower bounds on communication}
\label{sec:bounds}

In this section we focus on the parallel CALU\_PRRP algorithm based on
a binary reduction tree, and we show that it minimizes the
communication between different processors of a parallel computer.
For this, we use known lower bounds on the communication performed
during the LU factorization of a dense matrix of size $n \times n$,
which are
\begin{eqnarray}
\label{eq:lowerbound1}
\text{\# words\_moved} &=& \Omega\lt( \frac{n^3}{\sqrt M} \rt),  \\
\label{eq:lowerbound2}
\text{\# messages} &=& \Omega\lt( \frac{n^3}{M^{3\over2}} \rt), 
\end{eqnarray}
where \textit{\# words\_moved} refers to the volume of communication,
\textit{\# messages} refers to the number of messages exchanged, and
$M$ refers to the size of the memory (the fast memory in the case of a
sequential algorithm, or the memory per processor in the case of a
parallel algorithm).  These lower bounds were first introduced for
dense matrix multiplication \cite{hong1981io},
\cite{irony2004communication}, generalized later to LU factorization
\cite{demmel07:_tall_skinn_qr}, and then to almost all direct linear
algebra \cite{balard11:_min_com_linear_alg}.
\begin{comment}
with 2D layout. We first recall communication lower bounds for dense
matrix multiplication which was first introduced by Hong and Kung
\cite{hong1981io} in the sequential case (for not Strassen like
algorithms) and then generalized to the parallel case by Irony et
al. in \cite{irony2004communication}. Then Demmel et al. prove in
\cite{balard11:_min_com_linear_alg} that these lower bounds apply for
all direct linear algebra and so for the LU decomposition.

Let M be the local/ fast memory size, general lower bounds for the
$\#$ messages and the $\#$ words moved are:
\end{comment}
Note that these lower bounds apply to algorithms based on orthogonal
transformations under certain conditions
\cite{balard11:_min_com_linear_alg}.  However, this is not relevant to
our case, since CALU\_PRRP uses orthogonal transformations only to
select pivot rows, while the update of the trailing matrix is still
performed as in the classic LU factorization algorithm.  Hence the lower
bounds from equations \eqref{eq:lowerbound1} and, \eqref{eq:lowerbound2}
are valid for {\algB}.  

We estimate now the cost of computing in parallel the {\algB}
factorization of a matrix $A$ of size $m \times n$.  The matrix is
distributed on a grid of $P = P_r \times P_c$ processors using a
two-dimensional (2D) block cyclic layout. We use the following performance model. 
Let $\gamma$ be the cost
of performing a floating point operation, and let $\alpha + \beta w$
be the cost of sending a message of size $w$ words, where $\alpha$ is
the latency cost and $\beta$ is the inverse of the bandwidth. Then,
the total running time of an algorithm is estimated to be
\begin{eqnarray*}
\alpha \cdot (\text{\# messages}) + \beta \cdot (\text{\# words\_moved}) + \gamma \cdot (\text{\# flops} ),
\end{eqnarray*}	
where \#messages, \#words\_moved, and \#flops are counted along the
critical path of the algorithm.

Table \ref{tbl:CALUPRRP:par:model} displays the performance of
parallel {\algB} (a detailed estimation of the counts is presented in
Appendix D).  It also recalls the performance of two existing
algorithms, the PDGETRF routine from ScaLAPACK which implements GEPP,
and the CALU factorization.  All three algorithms have the same volume
of communication, since it is known that PDGETRF already minimizes the
volume of communication.  However, the number of messages of both
{\algB} and CALU is smaller by a factor of the order of $b$ than the
number of messages of PDGETRF.  This improvement is achieved thanks to
tournament pivoting.  In fact, partial pivoting, as used in the
routine PDGETRF, leads to an $O(n \log{P})$ number of messages, and
because of this, GEPP cannot minimize the number of messages.

Compared to CALU, {\algB} sends a small factor of less messages
(depending on $P_r$ and $P_c$) and performs $\frac{1}{P_r} \left(2 mn - n^2
   \right)  b +\frac{n b^2}{3}(5 \log_2{P_r} +1)$ more flops 
(which represents a lower order term).
This is because {\algB} uses the strong RRQR factorization at every node
of the reduction tree of every panel factorization, while CALU uses
GEPP.
\begin{comment}
The algorithm is performed in two steps.  The first is the computation
of the block LU factorization, the second is the computation of the
full LU factorization by performing additional GEPP on the $b\times b$
diagonal blocks.  The first step is also performing in two steps, a
preprocessing step to select the b pivot rows used for the panel
factorization, using Strong RRQR factorizations through the reduction
tree, then the application of the computed permutation matrix, the QR
factorization of the transpose of the current panel without pivoting
and finally the rank-b updates on the trailing matrix.
\end{comment}

\begin{table}[h]
\small \centering
\caption{Performance estimation of parallel (binary tree based) {\algB},
  parallel CALU, and PDGETRF routine when factoring an $m \times n$
  matrix, $m \geq n$.  The input matrix is distributed using a 2D
  block cyclic layout on a $P_r \times P_c$ grid of processors. 
Some lower order terms are omitted.}
\label{tbl:CALUPRRP:par:model}
\begin{tabular}{|l|l|l|l|}
\hline
  & Parallel CALU\_PRRP \\ \hline
\# messages & $ \frac{3n}{b} \log_2 P_r 
            +  \frac{2n}{b} \log_2{P_c}  $
           \\ 
\# words    & $    \left( nb + \frac{3}{2}\frac{ n^2}{ P_c} \right) \log_2 P_r +
        \frac{1}{P_r}\left(mn - \frac{n^2}{2} \right) \log_2{P_c}  $ 
           \\ 
\# flops    & $ \frac{1}{P} \left( mn^2 - \frac{n^3}{3} \right) + \frac{2}{P_r} \left(
2 mn - n^2   \right)  b + \frac{n^2b}{2 P_c}+   \frac{10n b^2}{3} \log_2{P_r}   $ 

\\ \hline

            & Parallel CALU  \\ \hline
\# messages & $ \frac{3n}{b} \log_2 P_r 
            +  \frac{3 n}{b} \log_2{P_c}  $
           \\ 
\# words    & $  \left( nb + \frac{3 n^2}{2 P_c} \right) \log_2 P_r +
        \frac{1}{P_r}\left(mn - \frac{n^2}{2} \right) \log_2{P_c}  $ 
           \\ 
\# flops    & $\frac{1}{P} \left( mn^2 - \frac{n^3}{3} \right) + \frac{1}{P_r} \left(
2 mn - n^2   \right)  b + \frac{n^2b}{2 P_c} + 
 \frac{n b^2}{3} (5 \log_2{P_r} -1) 
$ 
\\ \hline
         & PDGETRF \\ \hline
\# messages & $2n \left( 1 + \frac{2}{b} \right) \log_2{P_r} +
            \frac{3n}{b} \log_2{P_c} $ \\
\# words & $\left( \frac{nb}{2} +  \frac{3 n^2}{2P_c} \right)
        \log_2{P_r} + \log_2{P_c} \frac{1}{P_r}\left(mn -
        \frac{n^2}{2} \right)$ \\
\# flops   & $ \frac{1}{P} \left(m n^2 - \frac{n^3}{3} \right) +  \frac{1}{P_r} \left( mn - \frac{n^2}{2} \right) b
 + \frac{n^2 b}{2 P_c} $ \\
   \hline
\end{tabular}
\end{table}

Despite this additional communication cost, we show now that
CALU\_PRRP is optimal in terms of communication.  We choose optimal
values of the parameters $P_r, P_c$, and $b$, as used in CAQR
\cite{demmel07:_tall_skinn_qr} and CALU \cite{Grigori:EECS-2010-29},
that is,
\begin{equation*}
P_r =  \sqrt{\frac{m P}{n}}
\; \; , \;
P_c =  \sqrt{\frac{n P}{m}}
\; \; \text{and} \;
b   = \frac{1}{4}\log^{-2} \left( \sqrt{\frac{mP}{n}} \right) \cdot \sqrt{\frac{m n}{P}}=\log^{-2} \left(\frac{mP}{n} \right)\cdot \sqrt{\frac{m n}{P}}.
\end{equation*}
For a square matrix of size $n \times n$, the optimal parameters are,
\begin{equation*}
P_r =  \sqrt{P}
\; \; , \;
P_c =  \sqrt{P}
\; \; \text{and} \;
b   = \frac{1}{4}\log^{-2} \left( \sqrt{P} \right) \cdot \frac{n}{\sqrt{P}} = \log^{-2} \left( P \right) \cdot \frac{n}{\sqrt{P}}.
\end{equation*}

Table \ref{tbl:CALUPRRP:par:optimal} presents the performance
estimation of parallel {\algB} and parallel CALU when using the
optimal layout.  It also recalls the lower bounds on communication
from equations \eqref{eq:lowerbound1} and \eqref{eq:lowerbound2} when
the size of the memory per processor is on the order of $n^2/P$.  Both
{\algB} and CALU attain the lower bounds on the number of words and on
the number of messages, modulo polylogarithmic factors.  Note that the
optimal layout allows to reduce communication, while keeping the
number of extra floating point operations performed due to tournament
pivoting as a lower order term.  While in this section we focused on
minimizing communication between the processors of a parallel
computer, it is straightforward to show that the usage of a flat tree
during tournament pivoting allows {\algB} to minimize communication
between different levels of the memory hierarchy of a sequential
computer.

\begin{table}[h]
\small \centering
\caption{Performance estimation of parallel (binary tree based)
  CALU\_PRRP and CALU with an optimal layout.  The matrix factored is
  of size $n \times n$.  Some lower-order terms are omitted.}
\label{tbl:CALUPRRP:par:optimal}
\begin{tabular}{|l | l | l| }
 \hline
            & Parallel CALU\_PRRP with optimal layout & Lower bound \\ \hline

\# messages & ${5\over 2} \sqrt{P} \log^3{P} $ &  $\Omega ( \sqrt{ P } )$   \\ 
\# words    &  $  \frac{n^2}{ \sqrt{ P} }  \left( {1\over 2} \log^{-1}{P} + 
\log{P} \right)$
 &  $\Omega ( \frac{ n^2}{\sqrt{ P}}  )$ \\ 
\# flops   & $\frac{1}{P} \frac{2 n^3}{3}  +
\frac{5 n^3}{ 2 P \log^2  P } + \frac{5 n^3}{ 3 P \log^3  P } $  & $\frac{1}{P} \frac{2 n^3}{3}$  \\ \hline
  & Parallel CALU with optimal layout & Lower bound \\ \hline
  
\# messages & $3\sqrt{P} \log^3{P} $ & 
            $\Omega ( \sqrt{ P } )$ \\ 
\# words    & $  \frac{n^2}{ \sqrt{ P} }  \left( {1\over 2} \log^{-1}{P} + 
\log{P} \right)$
 &  $\Omega ( \frac{ n^2}{\sqrt{ P}}  )$ \\ 
\# flops    & $\frac{1}{P} \frac{2 n^3}{3}  +
\frac{3 n^3}{2 P \log^2 P } + \frac{5 n^3}{ 6 P \log^3 P } $ 
&  $\frac{1}{P} \frac{2 n^3}{3}$ \\ \hline
\end{tabular}
\end{table}

\section{Less stable factorizations that can also minimize communication}
\label{sec:related_work}

In this section, we present briefly two alternative algorithms that
are based on panel strong RRQR pivoting and that are conceived such
that they can minimize communication.  But we will see that they can
be unstable in practice.  These algorithms are also based on block
algorithms, that factor the input matrix by traversing panels of size
$b$.  The main difference between them and {\algB} is the panel
factorization, which is performed only once in the alternative
algorithms.

We present first a parallel alternative algorithm, which we refer to
as block parallel {\algA}.  At each step of the block factorization, the
panel is partitioned into $P$ block-rows $[A_0; A_1; \ldots;
  A_{P-1}]$.  The blocks below the diagonal $b \times b$ block of the
current panel are eliminated by performing a binary tree of strong
RRQR factorizations.  At the leaves of the tree, the elements below
the diagonal block of each block $A_i$ are eliminated using strong
RRQR.  The elimination of each such block row is followed by the
update of the corresponding block row of the trailing matrix.
\begin{comment}
That is, we select b pivot rows in each block performing the strong
rank revealing QR on the transpose of the current block, apply the
transpose of the computed permutation on the block of size ${m\over
  P}\times n$ and then update the trailing matrix of the block, which
is of size ${m\over P}-b\times n-b$.  
\end{comment}
The algorithm continues by performing the strong RRQR factorization of
pairs of $b \times b$ blocks stacked atop one another, until all the
blocks below the diagonal block are eliminated and the corresponding
trailing matrices are updated.  The algebra of the block parallel
{\algA} algorithm is detailed in Appendix E, while in figure \ref{fig:blockparallel_luprrp}
we illustrate one step of the factorization by using an arrow
notation, where the function $g(A_{ij})$ computes a strong RRQR on the matrix
$ A_{ij}^T$ and updates the trailing matrix in the same step. 

\begin{comment}
As for {\algB}, the strong RRQR to select pivot rows, but we update the
corresponding trailing matrix of each node during the same step. That
means at each node of the reduction tree we perform LU\_PRRP
factorization, and not select pivot rows during the preprocessing step
and then apply the permutation matrix and update the whole trailing
matrix as in CALU\_PRRP. For these algorithms, we do not perform the
QR factorization of the whole panel at the end as for CALU\_PRRP,
since the update and the computation of the L factor are done as we go
through the reduction tree.
\end{comment}

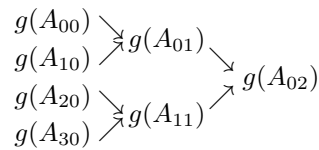
\begin{figure}[!htbp]
\begin{center}
\setlength{\unitlength}{.5cm}
\begin{picture}(7,4)
\put(0.5,0.5){$g(A_{30})$}
\put(0.5,1.5){$g(A_{20})$}
\put(0.5,2.5){$g(A_{10})$}
\put(0.5,3.5){$g(A_{00})$}

\put(2.7,0.65){$\nearrow$}
\put(2.7,1.35){$\searrow$}
\put(2.7,2.65){$\nearrow$}
\put(2.7,3.35){$\searrow$}

\put(3.5,1.0){$g(A_{11})$}
\put(3.5,3.0){$g(A_{01})$}

\put(5.6,1.5){$\nearrow$}
\put(5.6,2.5){$\searrow$}

\put(6.5,2.0){$g(A_{02})$}

\end{picture}

\caption{Block parallel LU\_PRRP\label{fig:blockparallel_luprrp}}
\end{center}
\end{figure}
\begin{figure}[!htbp]
\begin{center}
\setlength{\unitlength}{.5cm}
\begin{picture}(7,4)
\put(0.5,0.5){$A_{30}$}
\put(0.5,1.5){$A_{20}$}
\put(0.5,2.5){$A_{10}$}
\put(0.5,3.5){$A_{00}$}

\put(2.0,1.0){\vector(4,1){9}}
\put(2.0,1.75){\vector(4,1){6}}
\put(2.0,2.75){\vector(4,1){3}}
\put(1.8,3.75){\vector(1,0){.5}}

\put(2.2,3.5){$g(A_{00})$}
\put(4.4,3.75){\vector(1,0){.5}}
\put(5.0,3.5){$g(A_{01})$}
\put(7.2,3.75){\vector(1,0){.6}}
\put(8.0,3.5){$g(A_{02})$}
\put(10.2,3.75){\vector(1,0){.6}}
\put(11.0,3.5){$g(A_{03})$}
\end{picture}
\caption{Block pairwise LU\_PRRP \label{fig:blockpairwise_luprrp}}
\end{center}

\end{figure}
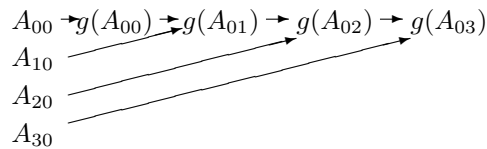

A sequential version of the algorithm is based on the usage of a flat
tree, and we refer to this algorithm as block pairwise {\algA}.  Using
the arrow notation, the figure \ref{fig:blockpairwise_luprrp} illustrates 
the elimination of one panel. 

The block parallel LU\_PRRP and the block pairwise LU\_PRRP algorithms
have similarities with the block parallel pivoting and the block
pairwise pivoting algorithms.  These two latter algorithms were shown
to be potentially unstable in~\cite{Grigori:EECS-2010-29}.  There is a
main difference between all these alternative algorithms and
algorithms that compute a classic LU factorization as GEPP, {\algA}, and
their communication avoiding variants.  The alternative algorithms
compute a factorization in the form of a product of lower triangular
factors and an upper triangular factor.  And the elimination of each
column leads to a rank update of the trailing matrix larger than one.
It is thought in~\cite{trefethen90:_averag_gauss} that the rank-1
update property of algorithms that compute an LU factorization
inhibits potential element growth during the factorization, while a
large rank update might lead to an unstable factorization.

Note however that at each step of the factorization, block parallel
and block pairwise {\algA} use at each level of the reduction tree
original rows of the active matrix.  Block parallel pivoting and block
pairwise pivoting algorithms use $U$ factors previously computed to achieve the
factorization, and this could potentially lead to a faster propagation of
ill-conditioning.

\begin{comment}
For these algorithms, the obtained $L$ factor is bounded by $\tau$,
since it is obtained through the different levels of the reduction
tree where each sub-block is bounded.  

So these algorithms are more stable than block parallel and block
pairwise pivoting methods which do not guarantee the fact that the
elements of L are bounded, a property which was shown again
in~\cite{trefethen90:_averag_gauss} to be important for the stability
of the factorization.  
\end{comment}

\begin{figure}[!htbp]
\centering
\includegraphics[scale=0.6]{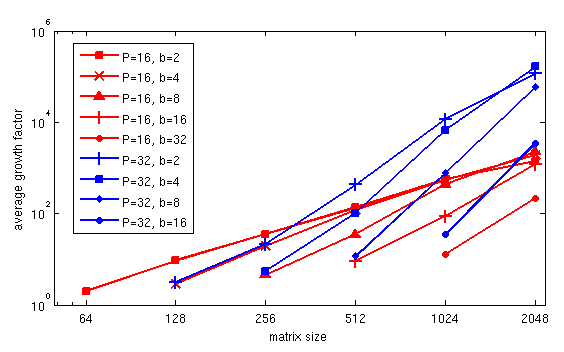}
\caption{\label{fig:bluprrp} Growth factor of block parallel LU\_PRRP
  for varying block size $b$ and number of processors $P$.}
\end{figure}

\begin{figure}[!htbp]\centering
\includegraphics[scale=0.43]{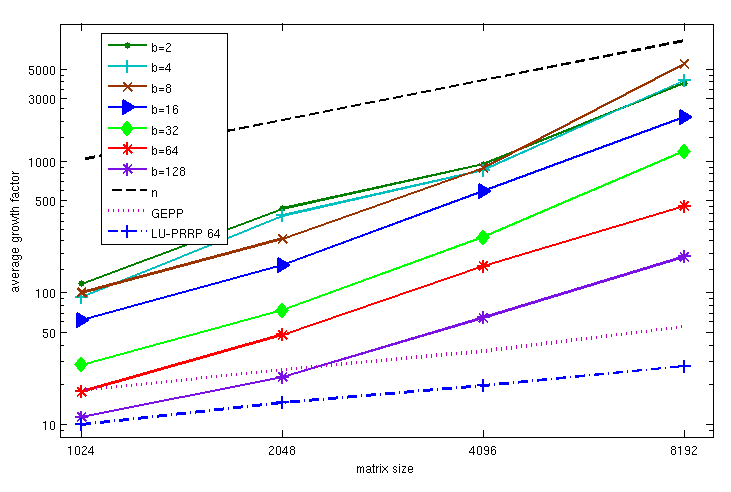}
\caption{\label{fig:fluprrp} Growth factor of block pairwise LU\_PRRP
  for varying matrix size and varying block size $b$.}
\end{figure}

The upper bound of the growth factor of both block parallel and block
pairwise {\algA} is $(1+\tau b)^{n\over b}$, since for every panel
factorization, a row is updated only once.  Hence they have the same
bounds as the {\algA} factorization, and smaller than that of the {\algB}
factorization.  Despite this, they are less stable than the {\algB}
factorization.  Figures \ref{fig:bluprrp} and \ref{fig:fluprrp}
display the growth factor of block parallel {\algA} and block pairwise
{\algA} for matrices following a normal distribution. In figure
\ref{fig:bluprrp}, the number of processors $P$ on which each panel is
partitioned is varying from 16 to 32, and the block size $b$ is varying
from $2$ to $16$.  The matrix size varies from $64$ to $2048$, but we
have observed the same behavior for matrices of size up to $8192$.
When the number of processors $P$ is equal to $1$, the block parallel
{\algA} corresponds to the LU\_PRRP factorization.  The results show
that there are values of $P$ and $b$ for which this method can be very
unstable.  For the sizes of matrices tested, when $b$ is chosen such
that the blocks at the leaves of the reduction tree have more than
$2b$ rows, the number of processors $P$ has an important impact, the
growth factor increases with increasing $P$, and the method is
unstable.

\begin{comment}
We can see that the size of the blocks used at the leaves of the
reduction tree has an important impact on the growth factor.  In fact,
when the number of rows in the leaves of the reduction tree is between
$b+1$ and $2b$, the method is very stable regardless of the number of
processors $P$ and the size of the panel $b$.  However when the number
of rows in the first level nodes of the reduction tree is bigger than
$2b$ the algorithm becomes unstable, the impact of the number of
processors $P$ becomes important and the growth factor increases with
increasing $P$.

This suggests that more investigations and further experiments are
necessary to understand the stability of binary tree LU\_PRRP
algorithm especially concerning the number of rows in the first level
nodes of the reduction tree.
\end{comment}

In Figure \ref{fig:fluprrp}, the matrix size varies from $1024$ to
$8192$.  For a given matrix size, the growth factor increases with
decreasing the size of the panel $b$, as one could expect.  We note
that the growth factor of block pairwise {\algA} is larger than that
obtained with the CALU\_PRRP factorization based on a flat tree scheme
presented in Table \ref{Tab_fcarrqr_randn1lu}.  But
it stays smaller than the size of the matrix $n$ for different panel
sizes.  Hence this method is more stable than block parallel {\algA}.
Further investigation is required to conclude on the stability of
these methods.
  
\newpage

\section{Conclusions}
\label{sec:conclusion}

This paper introduces {\algA}, an LU factorization algorithm based on
panel rank revealing pivoting.  This algorithm is more stable than
GEPP in terms of worst case growth factor.  It is also very stable in
practice for various classes of matrices, including pathological cases
on which GEPP fails.

Its communication avoiding version, {\algB}, is also more stable in
terms of worst case growth factor than CALU, the communication
avoiding version of GEPP.  More importantly, there are cases of
interest for which the upper bound of the growth factor of {\algB} is
smaller than that of GEPP for several cases of interest.  Extensive
experiments show that {\algB} is very stable in practice and leads to
results of the same order of magnitude as GEPP, sometimes even better.

Our future work focuses on two main directions.  The first direction
investigates the design of a communication avoiding algorithm that has
smaller bounds on the growth factor than that of GEPP in general.
\begin{comment}
the design of an hybrid CALU algorithm that combine the advantages of
GEPP (few flops) and those of Strong RRQR (stability thanks to the
bound ensured), that is to perform GEPP and some additional swaps on
the L factor as reduction operation.  The goal of this direction is to
design a communication avoiding algorithm that is more stable than
GEPP in terms of worst case pivot growth.
\end{comment}
The second direction focuses on estimating the performance of {\algB} on
parallel machines based on multicore processors, and comparing it with
the performance of CALU.

\bibliographystyle{siam}
{
\bibliography{carrqrlu}
}

%\pagebreak
\newpage
\section*{Appendix A}
\label{sec:appendixA}
We describe briefly strong RRQR introduced by M. Gu and S. Eisenstat
in~\cite{gu96:_strong_rank_revealing_qr}.  This factorization will be
used in our new LU decomposition algorithm, which aims to obtain an 
upper bound of the growth factor smaller than GEPP (see
Section~\ref{sec:rrqrlu_algebra}.)  Consider a given threshold $\tau>
1$ and an $h \times p$ matrix $B$ with $ p > h$, a Strong RRQR
factorization on a matrix $B$ gives (with an empty $(2,2)$ block)
\begin{eqnarray*}
 B^{T} \Pi = QR = Q 
\left[ 
\begin{array}{cc} 
 R_{11} &R_{12} \\
\end{array}
\right], 
\end{eqnarray*}
where $\norm{R_{11}^{-1}R_{12}}_{max} \leq \tau$, with
$\norm{\ .\ }_{max}$ being the biggest 
entry of a given matrix in absolute value. This
factorization can be computed by a classical QR
factorization with column pivoting followed by a limited
number of additional swaps and QR updates if necessary. \\

\begin{algorithm}[h!]
\caption{Strong RRQR}
\label{srrqr}
\begin{algorithmic}[1]
\State{Compute $ B^{T} \Pi = QR $ using the classical RRQR
  with column pivoting }
\While {there exist i and j such that $\abs{{(
      R_{11}^{-1}R_{12})}_{ij} }  > \tau$ }
\State{Set $\Pi = \Pi \Pi_{ij} $ and compute the QR factorization of $R$ $\Pi_{ij} $ (QR updates)}       
\EndWhile
\Ensure{$ B^{T} \Pi = QR $  with $\norm{R_{11}^{-1}R_{12}}_{max} \leq \tau$}
  \end{algorithmic}
\end{algorithm}
The while loop in Algorithm~\ref{srrqr} interchanges any
pairs of columns that can increase $\abs{det(R_{11})}$ by at
least a factor $\tau$. At most $O(\log_{\tau} n)$ such
interchanges are necessary before Algorithm~\ref{srrqr}
finds a strong RRQR factorization. The QR factorization of
$B^{T} \Pi$ can be computed numerically via efficient and
numerically stable QR updating
procedures. See~\cite{gu96:_strong_rank_revealing_qr} for
details.

\section*{Appendix B}
\label{sec:appendixB}
We present experimental results for the {\algA} factorization, the binary tree based {\algB}, and the flat tree based {\algB}.
We show results obtained for the LU decomposition and the linear solver. Tables \ref{TabRRQRLU_randn1lu}, \ref{Tab_fcarrqr_randn1lu}, 
and \ref{Tab_bcarrqr_randn1lu} display the results obtained for random matrices. They show the growth factor, 
the norm of the factor L and U and their inverses, and the relative error of the decomposition.\\
Tables \ref{TabGEPP_spec1}, \ref{TabLUPRRP_spec1}, \ref{TabFcaluprrp_spec}, and \ref{TabBcaluprrp_spec} display the results obtained for the special matrices presented 
in Table \ref{Tab_allTestMatrix}. The size of the tested matrices is n = 4096. For {\algA} and  flat tree based {\algB}, the size of the panel is $b = 8$.
For binary tree based {\algB} we use $P = 64$ and $b=8$, this means that the size of the matrices used at the leaves of the reduction tree is $64 \times 8$.\\
Tables \ref{bcaluprrp_gepp} and \ref{fcaluprrp_gepp} present results for the linear solver using binary tree based and flat tree based {\algB}, together with CALU and GEPP for
comparison.\\
The tables are presented in the following order.
\begin{itemize}
\item Table \ref{TabRRQRLU_randn1lu}: Stability of the LU decomposition for {\algA} and GEPP on random matrices.
\item Table \ref{TabGEPP_spec1}: Stability of the LU decomposition for GEPP on special matrices.
\item Table \ref{TabLUPRRP_spec1}: Stability of the LU decomposition for {\algA} on special matrices.
\item Table \ref{bcaluprrp_gepp}: Stability of the linear solver using binary tree based {\algB}, binary tree based CALU, and GEPP.
\item Table \ref{fcaluprrp_gepp}: Stability of the linear solver using flat tree based {\algB}, flat tree based CALU, and GEPP.
\item Table \ref{Tab_fcarrqr_randn1lu}: Stability of the LU decomposition for flat tree based {\algB} and GEPP on random matrices.
\item Table \ref{TabFcaluprrp_spec}: Stability of the LU decomposition for flat tree based {\algB} on special matrices.
\item Table \ref{Tab_bcarrqr_randn1lu}: Stability of the LU decomposition for binary tree based {\algB} and GEPP on random matrices.
\item Table \ref{TabBcaluprrp_spec}: Stability of the LU decomposition for binary tree based {\algB} on special matrices.
\end{itemize}

{\small
\begin{center}
\begin{longtable}{clp{3.9in}}
\caption{Special matrices in our test set.
\label{Tab_allTestMatrix} }
\\ \hline

 {\sf No. }      & {\sf Matrix}     & {\sf Remarks }  \\
 1   &   hadamard    &   Hadamard matrix, hadamard($n$), where $n$, $n/12$, or $n/20$ is power of 2. \\
 2   &   house   &   Householder matrix,  $A=\mathrm{eye}(n)-\beta*v*v'$, where  $[v,\beta,s]$ = gallery('house', randn($n$, 1)). \\
 3   &   parter  &   Parter matrix, a Toeplitz matrix with most of singular values near $\pi$. gallery('parter', $n$), or $A(i,j) = 1/(i-j+0.5)$.    \\
 4   &   ris     &   Ris matrix, matrix with elements $A(i,j) =0.5/(n-i-j+1.5)$. The eigenvalues cluster around $ -\pi/2$ and $\pi/2$. gallery('ris', $n$).  \\
 5   &   kms     &   Kac-Murdock-Szego Toeplitz matrix. Its inverse is tridiagonal.  gallery('kms', $n$) or gallery('kms', $n$, rand).  \\
 6   &   toeppen &   Pentadiagonal Toeplitz matrix (sparse). \\
\hline
\end{longtable}
\end{center}
}

{\small
 
\begin{landscape}\centering

\begin{table}[!htbp]\centering
\caption{Stability of the LU decomposition for {\algA} and GEPP on random matrices.}
\label{TabRRQRLU_randn1lu}
\begin{tabular}{|c|c||c|c|c|c|c|c|c|c|c|c|}
\hline
\multicolumn{9}{|c|}{{\algA} }\\
\hline {\sf n} & {\sf b} & $g_W$ &$||U||_1$ & $||U^{-1}||_1$ & $\frac{||PA-LU||_F}{||A||_F}$&HPL1&HPL2&HPL3  \\\hline
\multirow{5}{*}{8192}
        &    128  &2.38E+01  &9.94E+04    & 3.34E+02  &5.26E-14     &4.13E-02     &2.35E-02    &4.68E-03      \\ \cline{2-9}
        &    64  & 2.76E+01 &1.07E+05    &1.06E+03   &5.24E-14     &5.38E-02     &2.08E-02    &4.36E-03      \\ \cline{2-9}
        &    32  & 3.15E+01   & 1.20E+05 &3.38E+02 &5.14E-14 &3.95E-02 & 2.27E-02&4.54E-03\\ \cline{2-9}
        &   16   & 3.63E+01 &1.32E+05 & 3.20E+02&4.91E-14 &3.29E-02 &1.89E-02&3.79E-03\\ \cline{2-9}
        &    8    & 3.94E+01 &1.41E+05 &2.65E+02 &4.84E-14 &2.95E-02 & 1.69E-02 & 3.37E-03\\ \hline
\multirow{5}{*}{4096}
        &    128  &1.75E+01  &3.74E+04    & 4.34E+02  &2.69E-14     &3.28E-02     &2.21E-02    &4.58E-03      \\
\hline

\end{tabular}
\end{table}
\end{landscape}

{\small
\begin{landscape}\centering

\begin{table}[!htbp]\centering
\caption{Stability of the LU decomposition for GEPP on special matrices.}
\label{TabGEPP_spec1}
\scalebox{0.95}{
\begin{tabular}{|c|l|c|c|c|c|c|c|c|c|c|c|c|}\hline

%\multicolumn{2}{|c|}{ \multirow{2}{*}{matrix } } & \multirow{2}{*}{cond(A,2) } & \multicolumn{10}{c|}{ GEPP } \\\cline{4-13}
%\multicolumn{2}{|c|}{} & & $g_W$  & $||L||_1$ & $||L^{-1}||_1$ & ${\max\atop{ij}} |U_{ij}|$ & ${\min\atop{kk}} |U_{kk}|$  & cond($U$,1) & $\frac{||PA-LU||_F}{||A||_F}$ & $\eta$ & $w_b$  & $N_{IR}$ \\ \hline

  & matrix & cond(A,2) & $g_W$  & $||L||_1$ & $||L^{-1}||_1$ & ${\max\atop{ij}} |U_{ij}|$ & ${\min\atop{kk}} |U_{kk}|$ & cond($U$,1) & $\frac{||PA-LU||_F}{||A||_F}$ & $\eta$ & $w_b$  & $N_{IR}$ \\
\hline

\multirow{26}{*}{ \begin{sideways} well-conditioned \end{sideways} }
    &   hadamard    &   1.0E+0  &   4.1E+3  &   4.1E+3  &   4.1E+3  &   4.1E+3  &   1.0E+0  &   5.3E+5  &   0.0E+0  &   3.3E-16 &   4.6E-15 &   2   \\ \cline{2-13}
    &   house   &   1.0E+0  &   5.1E+0  &   8.9E+2  &   2.6E+2  &   5.1E+0  &   5.7  &   1.4E+4  &   2.0E-15 &   5.6E-17 &   6.3E-15 &   3   \\ \cline{2-13}
    &   parter  &   4.8E+0  &   1.6E+0  &   4.8E+1  &   2.0E+0  &   3.1E+0  &   2.0E+0  &   2.3E+2  &   2.3E-15 &   8.3E-16 &   4.4E-15 &   3   \\ \cline{2-13}
    &   ris &   4.8E+0  &   1.6E+0  &   4.8E+1  &   2.0E+0  &   1.6E+0  &   1.0E+0  &   2.3E+2  &   2.3E-15 &   7.1E-16 &   4.7E-15 &   2   \\ \cline{2-13}
    &   kms &   9.1E+0  &   1.0E+0  &   2.0E+0  &   1.5E+0  &   1.0E+0  &   7.5E-1  &   3.0E+0  &   2.0E-16 &   1.1E-16 &   6.7E-16 &   1   \\ \cline{2-13}
    &   toeppen &   1.0E+1  &   1.1E+0  &   2.1E+0  &   9.0E+0  &   1.1E+1  &   1.0E+1  &   3.3E+1  &   1.1E-17 &   7.2E-17 &   3.0E-15 &   1   \\ \cline{2-13}
    &   condex  &   1.0E+2  &   1.0E+0  &   2.0E+0  &   5.6E+0  &   1.0E+2  &   1.0E+0  &   7.8E+2  &   1.8E-15 &   9.7E-16 &   6.8E-15 &   3   \\
\hline

\end{tabular}}
\end{table}

\end{landscape}

}
%Special matrices
{\small
\begin{landscape}\centering
\begin{table}[!htbp]\centering
\caption{Stability of the LU decomposition for {\algA} on special matrices.}
\label{TabLUPRRP_spec1}
\scalebox{0.95}{\begin{tabular}{|c|l|c|c|c|c|c|c|c|c|c|c|}\hline
  & matrix & cond(A,2) & $g_W$ & $||L||_1$ & $||L^{-1}||_1$  & $||U||_1$ & $||U^{-1}||_1$ & $\frac{||PA-LU||_F}{||A||_F}$ & $\eta$ & $w_b$  & $N_{IR}$ \\
\hline

\multirow{26}{*}{ \begin{sideways} well-conditioned \end{sideways} }
    &   hadamard    &   1.0E+0  &  5.1E+02 &  5.1E+02 &  2.5E+16&   5.6E+0& 6.0E+00&   5.2E-15 &  1.1E-15 &  9.0E-15 &  2 \\ \cline{2-12}
    &   house   &   1.0E+0  &  7.4E+00 &   8.3E+02 &  7.9E+02 &  3.7E+02 &  5.0E+01 &  5.8E-16 &  4.2E-17  & 4.8E-15  & (2) \\ \cline{2-12}
    &   parter  &   4.8E+0  &  1.5E+00  & 4.7E+01  & 3.7E+00 &  1.4E+01 &  3.6E+01 &  8.5E-16  & 6.8E-16 &  4.4E-15 &  3 \\ \cline{2-12}
    &   ris &   4.8E+0  &  1.5E+00  & 4.7E+01  & 1.1E+01 &  7.3E+00 &  7.2E+01  & 8.5E-16  & 7.6E-16 &  4.6E-15 &  2 \\ \cline{2-12}
    &   kms &   9.1E+0  &  1.0E+00  & 6.7E+00  & 4.0E+00 &  4.1E+00 &  1.2E+01 &  1.4E-16 &  6.8E-17  & 3.7E-16  & 1 \\ \cline{2-12}
    &   toeppen &   1.0E+1  &  1.3E+00  & 2.3E+00  & 3.4E+00  & 3.6E+01 &  1.0E+00  & 1.3E-16  & 8.3E-17 &  2.2E-15 &  1 \\ \cline{2-12}
    &   condex  &   1.0E+2  &  1.0E+00 &  1.9E+00  & 5.5E+00 &  5.9E+02  & 1.2E+00 &  6.5E-16 &  7.4E-16  & 5.1E-15 &  (2)\\
\hline

\end{tabular}}
\end{table}
\end{landscape}
}

{\small
\centering
\begin{table}[!htbp]\centering
\caption{Stability of the linear solver using binary tree based {\algB}, binary tree based CALU, and GEPP.}
\label{bcaluprrp_gepp}
\begin{tabular}{|c|c|c||c|c|c|c|c|c|}
\hline {\sf n} & {\sf P} & {\sf b} &  $\eta$ & $w_b$   & $N_{IR}$ &
{\sf HPL1} & {\sf HPL2} & {\sf HPL3}  \\  \hline

\multicolumn{9}{|c|}{Binary tree based \algB}\\ \hline

\multirow{9}{*}{8192}
        &   \multirow{2}{*}{256} &   32  & 7.5E-15  &4.4E-14   &2     &6.2E-02   &2.3E-02    &4.4E-03    \\ \cline{3-9}
        &       &   16  & 6.7E-15   & 4.1E-14  & 2    &3.4E-02   &2.2E-02    &4.6E-03   \\ \cline{2-9}

        &   \multirow{3}{*}{128} &   64  &  7.6E-15  &4.7E-14   &2     &2.6E-02   &2.5E-02    &4.9E-03   \\ \cline{3-9}
        &       &   32  &  7.5E-15 & 4.9E-14  & 2    &6.9E-02   &2.6E-02    &4.9E-03 \\ \cline{3-9}
        &       &   16  & 7.3E-15   & 5.1E-14  &  2   &  4.2E-02 &  2.7E-02  &5.7E-03   \\ \cline{2-9}

        &   \multirow{4}{*}{64}  &   128 & 7.6E-15  & 5.0E-14  &2     &2.8E-02   &2.6E-02    &5.2E-03   \\
\hline

\end{tabular}
\end{table}

}

{\small
\centering
\begin{table}[!htbp]\centering
\caption{Stability of the linear solver using flat tree based {\algB},
 flat tree based CALU, and GEPP.}
\label{fcaluprrp_gepp}
\begin{tabular}{|c|c|c||c|c|c|c|c|c|}
\hline {\sf n} & {\sf P} & {\sf b} &  $\eta$ & $w_b$   & $N_{IR}$ &
{\sf HPL1} & {\sf HPL2} & {\sf HPL3}  \\  \hline

\multicolumn{9}{|c|}{Flat tree based \algB}\\ \hline
 \multirow{4}{*}{8096}   
        &   -   &   8   & 6.2E-15  &3.8E-14   &1    &1.7E-02   &2.0E-02     &4.0E-03     \\ \cline{2-9}
        &   -   &   16  & 6.6E-15  &4.3E-14   &1    &3.6E-02   &2.2E-02     &4.3E-03 \\ \cline{2-9}
        &   -   &   32  & 7.2E-15  & 4.8E-14 &1    &6.5E-02   &2.5E-02     &5.3E-03 \\ \cline{2-9}
        &   -   &   64  & 7.3E-15  &5.1E-14   &1    &4.8E-02   &2.7E-02     &5.4E-03 \\ \hline

    \multirow{4}{*}{4096}  
        &   -   &   8   & 3.2E-15  &2.0E-14   &1    &6.7E-02   &2.1E-02     &4.7E-03 \\ \cline{2-9}
        &   -   &   16  & 3.6E-15  &2.3E-14  &1    &2.9E-02   &2.4E-02     &5.1E-03 \\
\hline

\end{tabular}
\end{table}
 
}

\begin{landscape}\centering
\begin{table}[!htbp]\centering
\caption{Stability of the LU decomposition for flat tree based {\algB} and GEPP on random matrices}
\label{Tab_fcarrqr_randn1lu}
\scalebox{0.95}{\begin{tabular}{|c|c||c|c|c|c|c|c|c|c|c|c|}
\hline
\multicolumn{11}{|c|}{ Flat Tree based {\algB}}\\
\hline {\sf n} & {\sf b} & $g_W$ &$||L||_1$ & $||L^{-1}||_1$  &
$||U||_1$ & $||U^{-1}||_1$ & $\frac{||PA-LU||_F}{||A||_F}$&HPL1&HPL2&HPL3  \\\hline
\multirow{5}{*}{8192}
        &    128  &3.44E+01  &2.96E+03&  2.03E+03&   1.29E+05&  2.86E+03&   8.67E-14&   1.77E-01&  2.53E-02&  5.00E-03\\ \cline{2-11}
        &    64  &4.42E+01&  3.30E+03&  2.28E+03&  1.43E+05&  6.89E+02&  8.56E-14&   3.62E-02&  2.97E-02&  6.40E-03\\ \cline{2-11}
        &    32  & 5.88E+01  &4.34E+03  & 2.50E+03   & 1.54E+05   & 3.79E+02  &8.15E-14  &3.40E-02  &2.96E-02  & 6.40E-03\\ \cline{2-11}
        &   16 & 6.05E+01 &4.54E+03  &2.50E+03  &1.53E+05  &4.57E+02  &6.59E-14  &5.66E-02  &2.27E-02  &4.60E-03 \\ \cline{2-11}
        &    8 & 5.25E+01 &3.24E+03 &2.63E+03 &1.60E+05 &1.78E+02 &5.75E-14 &1.92E-02 &1.93E-02 &4.10E-03 \\ \hline

\multirow{5}{*}{4096}
        &    128  &2.38E+01&1.60E+03&1.05E+03&4.75E+04& 5.11E+03&4.15E-14&1.58E+00&2.40E-02&5.16E-03\\
\hline

\end{tabular}}
\end{table}
\end{landscape}

%Special matrices
{\small
\begin{landscape}\centering
\begin{table}[!htbp]\centering
\caption{Stability of the LU decomposition for flat tree based {\algB} on special matrices.} \label{TabFcaluprrp_spec}
\scalebox{0.95}{\begin{tabular}{|c|l|c|c|c|c|c|c|c|c|c|c|}\hline
  & matrix & cond(A,2) & $g_W$ & $||L||_1$ & $||L^{-1}||_1$  & $||U||_1$ & $||U^{-1}||_1$ & $\frac{||PA-LU||_F}{||A||_F}$ & $\eta$ & $w_b$  & $N_{IR}$ \\
\hline

\multirow{26}{*}{ \begin{sideways} well-conditioned \end{sideways} }
    &   hadamard    &   1.0E+0  & 5.1E+02 &  5.1E+02 &  1.4E+02 &  5.8E+04 &  5.3E+00 &  5.2E-15  & 9.5E-16  & 7.9E-15  & 1 \\ \cline{2-12}
    &   house   &   1.0E+0  & 6.5E+00 &  9.0E+02 &  3.3E+02  & 3.3E+02  & 5.2E+01 &  6.6E-16 &  4.8E-17 &  5.3E-15 &  1 \\ \cline{2-12}
    &   parter  &   4.8E+0  & 1.5E+00  & 4.7E+01 &  3.7E+00 &  1.4E+01  & 3.6E+01 &  8.6E-16 &  6.8E-16 &  4.4E-15 &  1 \\ \cline{2-12}
    &   ris &   4.8E+0  &  1.5E+00 &  4.7E+01  & 3.7E+00 &  7.3E+00 &  7.2E+01 &  8.5E-16 &  6.8E-16 &  4.1-15  & 1 \\ \cline{2-12}
    &   kms &   9.1E+0  & 1.0E+00 &  5.4E+00  & 3.2E+00  & 3.9E+00  & 9.8E+00 &  1.13E-16 &  7.0E-17 &  4.2E-16 &  1 \\ \cline{2-12}
    &   toeppen &   1.0E+1  & 1.3E+00 &  2.3E+00  & 3.4E+00  & 3.6E+01 &  1.0E+00 &  1.0E-16 &  8.3E-17  & 2.7E-15 &  1 \\ \cline{2-12}
    &   condex  &   1.0E+2  & 1.0E+00  & 1.9E+00  & 5.5E+00 &  5.9E+02 &  1.2E+00 &  6.5E-16 &  7.4E-16 &  5.2E-15  & 1 \\
\hline

\end{tabular}}
\end{table}
\end{landscape}
}

\begin{landscape}\centering
\begin{table}[!htbp]\centering
\caption{Stability of the LU decomposition for binary tree based  {\algB} and GEPP on random matrices.}
\label{Tab_bcarrqr_randn1lu}
\begin{tabular}{|c|c|c||c|c|c|c|c|c|c|c|c|c|}
\hline
\multicolumn{12}{|c|}{ Binary Tree based {\algB} }\\
\hline {\sf n} & {\sf P} & {\sf b} & $g_W$ &$||L||_1$ & $||L^{-1}||_1$  &
$||U||_1$ & $||U^{-1}||_1$ & $\frac{||PA-LU||_F}{||A||_F}$&HPL1&HPL2&HPL3  \\\hline

\multirow{10}{*}{8192}
        &   256 &   16  &4.29E+01  & 4.41E+03   &  2.68E+03 & 1.62E+05 & 2.15E+02  &  7.39E-14  &  3.44E-02  & 2.20E-02 & 4.69E-03  \\ \cline{2-12}
     
        &   \multirow{2}{*}{128}
                &   32  & 4.58E+01 & 3.47E+03   & 2.58E+03  & 1.56E+05 & 2.96E+02 &  8.67E-14  &  6.92E-02  & 2.63E-02&  4.98E-03 \\ \cline{3-12}
        &       &   16 &4.97E+01  &   4.05E+03 &  2.73E+03 & 1.63E+05 & 2.16E+02  &  7.59E-14  &  4.22E-02  & 2.71E-02& 5.76E-03 \\ \cline{2-12}
        &   \multirow{4}{*}{64}  
                &   64  &4.77E+01  &   3.23E+03 & 2.35E+03  & 1.45E+05 & 3.22E+02 & 9.14E-14   &  6.37E-02  & 2.85E-02 & 5.91E-03  \\ \cline{3-12}
        &        &   32  & 4.86E+01 &  4.16E+03  &  2.64E+03 & 1.58E+05 & 5.15E+02 &  9.04E-14  &  7.07E-02  & 2.69E-02&  5.09E-03\\ \cline{3-12}
        &       &   16  & 5.37E+01 &  3.83E+03  & 2.67E+03  & 1.64E+05 & 2.56E+02 & 7.30E-14   &   4.21E-02 &2.70E-02 &  5.75E-03\\
\hline

\end{tabular}
\end{table}

\end{landscape}

%Special matrices
{\small
\begin{landscape}\centering
\begin{table}[!htbp]\centering
\caption{Stability of the LU decomposition for binary tree based {\algB} on special matrices.} \label{TabBcaluprrp_spec}
\scalebox{0.95}{\begin{tabular}{|c|l|c|c|c|c|c|c|c|c|c|c|}\hline
  & matrix & cond(A,2) & $g_W$ & $||L||_1$ & $||L^{-1}||_1$  & $||U||_1$ & $||U^{-1}||_1$ & $\frac{||PA-LU||_F}{||A||_F}$ & $\eta$ & $w_b$  & $N_{IR}$ \\
\hline

\multirow{26}{*}{ \begin{sideways} well-conditioned \end{sideways} }
    &   hadamard    &   1.0E+0  & 5.1E+02 &  5.1E+02 &  1.3E+02 &  4.9E+04  & 6.5E+00 &  5.3E-15 &  1.2E-15  & 8.7E-15 &  2 \\ \cline{2-12}
    &   house   &   1.0E+0  & 7.4E+00  & 8.3E+02  & 3.8E+02  & 3.7E+02 &  5.0E+01 &  5.8E-16  & 4.2E-17  & 4.8E-15 &  3\\ \cline{2-12}
    &   parter  &   4.8E+0  &1.5E+00  & 4.7E+01  & 3.7E+00 &  1.4E+01 &  3.6E+01 &  8.6E-16  & 6.9E-16 &  4.4E-15 &  3 \\ \cline{2-12}
    &   ris &   4.8E+0  &  1.5E+00 &  4.7E+01 &  3.7E+00  & 7.3E+00  & 7.2E+01  & 8.5E-16 &  6.2E-16&   4.4E-15  & 2 \\ \cline{2-12}
    &   kms &   9.1E+0  & 1.0E+00  & 5.1E+00  & 3.4E+00 &  4.3E+00 &  9.3E+00 &  1.1E-16 &  6.3E-17  & 3.6E-16 &  1\\ \cline{2-12}
    &   toeppen &   1.0E+1  & 1.3E+00 &  2.3E+00 &  3.4E+00 &  3.6E+01&   1.0E+00  & 1.0E-16  & 8.2E-17 &  4.0E-15  & 1 \\ \cline{2-12}
    &   condex  &   1.0E+2  & 1.0E+00 &  1.9E+00  & 5.5E+00 &  5.9E+02 &  1.2E+00 &  6.5E-16  & 7.4E-16 &  5.2E-15  & 2 \\
\hline

\end{tabular}}
\end{table}
\end{landscape}
}

\section*{Appendix C}
\label{sec:appendixC}
Here we summarize the floating-point operation counts for the {\algA}
algorithm performed on an input matrix A of size $m\times n$, where
 $m\geq n$ . We first focus on the step k of  the algorithm, that 
is we consider the $k^{th}$ panel of size $(m-(k-1)b) \times b$. 
We first perform a Strong RRQR on the transpose of the considered 
panel, then update the trailing matrix of size $ (m-kb) \times
(n-kb)$, finally we perform GEPP on the diagonal block of size $b
\times b$.
We assume that $m-(k-1)b \geq b+1$. The Strong RRQR performs nearly 
as many floating-point operations as the QR with column pivoting. 
Here we consider that is the same, since in practice performing QR 
with column pivoting is enough to obtain the bound $\tau$, and thus
 it is 
 \begin{eqnarray*}
Flops_{SRRQR, 1 block, step k} =  2(m-(k-1)b)b^2 -{ 2\over 3}  b^3.
\end{eqnarray*}
If we consider the update step (\ref{update}), then the flops count is
 \begin{eqnarray*}
Flops_{update, step k} = 2b (m-kb)(n-kb).
\end{eqnarray*}
For the additional GEPP on the diagonal block, the flops count is
\begin{eqnarray*}
Flops_{gepp, step k} = { 2\over 3}  b^3 + (n-kb)b^2 .
\end{eqnarray*}
Then the flops count for the step k is
\begin{eqnarray*}
Flops_{\algA, step k} = b^2(2m+n-3(k-1)b) -  b^3+ 2b (m-kb)(n-kb) .
\end{eqnarray*}
This gives us  an arithmetic operation count of
\begin{eqnarray*}
Flops_{\algA}(m,n,b) &=& \Sigma_{k=1}^{n\over b }  \left[ b^2(2m+n-2(k-1)b-kb) + 2b (m-kb)(n-kb) \right],\\
Flops_{\algA}(m,n,b) &=&  mn^2+2mnb+2nb^2-{1\over 2}n^2b -{1\over 3}n^3\\
&\sim& mn^2-{1\over 3}n^3+2mnb-{1\over 2}n^2b.
\end{eqnarray*}
Then for a square matrix of size $n\times n$, the flops count is
 \begin{eqnarray*}
Flops_{\algA}(n,n,b) = { 2\over 3} n^3+{ 3\over 2}n^2b+2nb^2  \sim
{ 2\over 3}  n^3 + { 3\over 2}n^2b.
\end{eqnarray*}
\section*{Appendix D}
\label{sec:appendixD}
Here we detail the performance model of the parallel version of the {\algB}
factorization performed on an input matrix A of size $m\times n$ where, $m\geq n$. 
We consider a 2D layout $P = P_r\times P_c$ . We first focus on the panel 
factorization for the block LU factorization, that is the
selection of the b pivot rows with the tournament pivoting strategy. This step
 is similar to CALU except that the reduction operator is Strong RRQR instead of GEPP, then 
for each panel the amount of communication is the same as for TSLU:
\begin{eqnarray*}
 \text{\# messages} &=& \log P_r\\
 \text{\# words} &=& b^2 \log P_r 
\end{eqnarray*} 
However the floating-point operations count is different. We 
consider as in Appendix C that Strong RRQR performs as many flops as QR 
with columns pivoting, then the panel factorization performs 
the QR factorization with columns pivoting on the transpose of the blocks of the panel and 
$\log P_r$ reduction steps:

 \begin{eqnarray*}
 \text{\# flops} &=& 2 {{m-(k-1)b}\over P_r} b^2 -{ 2\over 3}  b^3 + \log P_r(2(2b)b^2 - { 2\over 3}  b^3) \\
&=& 2 {{m-(k-1)b} \over P_r}  b^2+{10 \over 3}b^3\log P_r -{ 2\over 3}  b^3
\end{eqnarray*}
To perform the QR factorization without pivoting on the transpose of the panel, and the update of the trailing matrix: 
\begin{itemize}
\item broadcast the pivot information along the rows of the process grid.
\begin{eqnarray*}
\text{\# messages} &=& \log P_c\\
\text{\# words} &=& b \log P_c
\end{eqnarray*}
\item apply the pivot information to the original rows.
 \begin{eqnarray*}
\text{\# messages} &=& \log P_r\\
\text{\# words} &=& {{nb}\over P_c}  \log P_r 
 \end{eqnarray*}
\item Compute the block column L and broadcast it through blocks of columns
\begin{eqnarray*}
\text{\# messages} &=& \log P_c \\
\text{\# words} &=&  {{m-kb}\over P_r}b \log P_c\\
\text{\# flops} &=& 2 {{m-kb}\over P_r} b^2
\end{eqnarray*}
\item broadcast the upper block of the permuted matrix A through blocks of rows
 \begin{eqnarray*}
\text{\# messages} &=& \log P_r \\
\text{\# words} &=&  {{n-kb}\over P_c}b \log P_r
\end{eqnarray*}
\item perform a rank-b update of the trailing matrix 
 \begin{eqnarray*}
\text{\# flops} = 2b {{m-kb}\over P_r} {{n-kb}\over P_c}
\end{eqnarray*}
\end{itemize}
Thus to get the block LU factorization :
 \begin{eqnarray*}
\text{\# messages} &=& {{3n}\over b} \log P_r + {{2n}\over b} \log P_c\\
\text{\# words   } &=&  ({mn\over P_r}-{1\over 2}{n^2\over P_r}+n)\log P_c+(nb+{3\over 2}{n^2\over P_c})\log P_r\\
\text{\# flops   } &=& {1\over P}(mn^2-{1\over 3}n^3)+{2\over 3}nb^2(5\log P_r -1)+{b\over P_r}(4mn+2nb-2n^2)\\
&\sim& {1\over P}(mn^2-{1\over 3}n^3)+{2\over 3}nb^2(5\log P_r -1)+{4\over P_r}(mn-{n^2\over 2})b
 \end{eqnarray*}
Then to obtain the full LU factorization, for each $b \times b$ block, we 
perform the Gaussian elimination with partial pivoting and we update the
corresponding trailing matrix of size $b\times n-kb$. During this 
additional step, we first perform GEPP on the diagonal block, broadcast pivot rows
through column blocks (this broadcast can be done together with the previous broadcast of L, thus
there is no additional message to send), apply pivots, and finally compute U.
%\begin{eqnarray*}
%&\text{\# messages} = \log P_c
%\end{eqnarray*}
\begin{eqnarray*}
&\text{\# words} = n(1+{b\over 2}) \log P_c
\end{eqnarray*}
\begin{eqnarray*}
\text{\# flops} = \sum_{k=1}^{n\over b} \left[{2\over 3}b^3+ {{n-kb}\over P_c}b^2\right] = {2\over 3}nb^2+ {1\over 2}{{n^2b}\over P_c}
 \end{eqnarray*}
Finally the total count is :
\begin{eqnarray*}
\text{\# messages} &=&{{3n}\over b} \log P_r + {{2n}\over b} \log P_c \\
\text{\# words} &=& ({mn\over P_r}-{1\over 2}{n^2\over P_r}+{nb\over 2}+2n)\log P_c+(nb+{3\over 2}{n^2\over P_c})\log P_r\\
&\sim& ({mn\over P_r}-{1\over 2}{n^2\over P_r}) \log P_c+(nb+{3\over 2}{n^2\over P_c})\log P_r\\
\text{\# flops} &=&   {1\over P}(mn^2-{1\over 3}n^3)+{4\over P_r}(mn-{n^2\over 2})b + {{n^2b}\over{2 P_c }}+{10\over 3}nb^2\log P_r
 \end{eqnarray*}

 \section*{Appendix E}
 \label{sec:appendixE}
Here we detail the algebra of the block parallel LU-PRRP.
 At the first iteration, the matrix A has the following partition
 $$
 A = 
  \left[
 \begin{array}{ccccc} 	
 A_{11} & A_{12}\\
 A_{21} & A_{22} \\
 \end{array} 	
 \right]. 
 $$ 
 The block $A_{11}$ is of size $m/p \times b$, where p is the number of
 processors used and $m/p \ge b+1$, the block $A_{12}$ is of size 
$m/p \times n-b$, the block $A_{21}$ is of size $m-m/p \times b$, and the 
block $A_{22}$ is of size $m-m/p \times n-b$.\\

 To describe the algebra, we consider 4 processors and a block of size $b$.
 We first focus on the $b$ first columns,
 $$
 A(:,1:b) = 
  \left[
 \begin{array}{c} 	
 A_{0}\\
 A_{1} \\
 A_{2} \\
 A_{3} \\
 \end{array} 	
 \right]. 
 $$ 
 Each block $A_{i}$ is of size $m/p \times b$.  
 In the following we describe the different steps of the panel factorization. First, 
we perform Strong RRQR factorization on the transpose of each block  $A_{i}$ so we obtain : \\

 $$
 \begin{array}{c}  
 A_{0}^T \Pi_{00} = Q_{00} R_{00} \\
 A_{1}^T \Pi_{10} = Q_{10} R_{10} \\
 A_{2}^T \Pi_{20} = Q_{20} R_{20} \\
 A_{3}^T \Pi_{30} = Q_{30} R_{30} \\

 \end{array} 
 $$

 Each matrix $R_i$ is of size $b \times m/p $. Using
 MATLAB notations, we can write $R_i$ as following :\\

 \begin{center}
 $\begin{array}{cc} 
 \bar{R}_{i} = R_{i} (1:b,1:b)&\bar{\bar{R}}_{i} = R_{i} (1:b,b+1:m/p) \\
 \end{array} 
 $
 \end{center}	
 This step aims to eliminate the last $m/p - b$ rows of each block $A_{i}$. We define
 the matrix $ D_0 \Pi_0$ : 
 {\tiny
 $$ D_0 \Pi_0 =
 \left[
 \begin{array}{cccc} 
 \left[ 
 \begin{array}{cc} 
 I_b& \\
 -D_{00}&I_{m/p - b}  \\
 \end{array} 
 \right] & & &\\
  &\left[  \begin{array}{cc} 
 I_b&  \\
 -D_{10}&I_{m/p - b}  \\
 \end{array} 
 \right] & & \\
  & & \left[ \begin{array}{cc} 
 I_b& \\
 -D_{20}&I_{m/p - b}  \\
 \end{array} 
 \right]  &  \\
 & & & \left[ 
 \begin{array}{cc} 
 I_b& \\
 -D_{30}&I_{m/p - b}  \\
 \end{array} 
 \right] \\
 \end{array} 
 \right]  \times
 \left[
 \begin{array}{cccc} 
 \Pi_{00}^T& & &\\
 &\Pi_{10}^T& &\\
 & &\Pi_{20}^T&\\
 & & &\Pi_{30}^T\\
 \end{array} 
 \right],  
 $$
 } 

 where \\
 $$
 \begin{array}{c}  
 D_{00}=\bar{\bar{R}}_{00}^T  (\bar{R}_{00}^{-1})^T,\\
 D_{10}=\bar{\bar{R}}_{10}^T  (\bar{R}_{10}^{-1})^T,\\
 D_{20}=\bar{\bar{R}}_{20}^T  (\bar{R}_{20}^{-1})^T,\\
 D_{30}=\bar{\bar{R}}_{30}^T (\bar{R}_{30}^{-1})^T.\\
 \end{array} 
 $$
 Multiplying  $A(:,1:b)$ by  $D_0 \Pi_0 $ we obtain : \\

 $$
 D_0 \Pi_0 \times A(:,1:b) =
 \left[ 
 \begin{array}{c}
  (\Pi_{00}^T\times
 A_{0})(1:b,1:b)  \\
 0_{m/p - b}\\
 \\
  (\Pi_{10}^T\times
 A_{1})(1:b,1:b)   \\
 0_{m/p - b} \\
 \\
  (\Pi_{20}^T\times A_{2})(1:b,1:b) \\
 0_{m/p - b}\\
 \\
  (\Pi_{30}^T\times A_{3})(1:b,1:b)   \\
 0_{m/p - b}\\
 \\
 \end{array} 
 \right] =
 \left[ 
 \begin{array}{c}
 A_{01}  \\
 0_{m/p - b}\\
 \\
 A_{11}\\
 0_{m/p - b} \\
 \\
 A_{21}\\
 0_{m/p - b}\\
 \\
 A_{31}\\
 0_{m/p - b}\\
 \\
 \end{array} 
 \right].
 $$

 The second step corresponds to the second level of the reduction tree. 
We merge pairs of $b \times b$ blocks and as in the previous step we perform
 Strong RRQR factorization on the transpose of the $2b \times b$ blocks :\\
 $$
 \left[ 
 \begin{array}{c}
 A_{01}  \\
 A_{11}\\
 \end{array} 
 \right] 
 $$
 and \\
 $$
 \left[ 
 \begin{array}{c}
 A_{21}\\
 A_{31}\\
 \end{array} 
 \right].
 $$
 We obtain \\

 $$
 \begin{array}{c}  
 \left[ 
 \begin{array}{c}
 A_{01}  \\
 A_{11}\\
 \end{array} 
 \right] ^T \bar{\Pi}_{01} = Q_{01} R_{01}, \\
 \left[ 
 \begin{array}{c}
 A_{21}\\
 A_{31}\\
 \end{array} 
 \right] ^T \bar{\Pi}_{11} = Q_{11} R_{11}. \\
 \end{array} 
 $$
 We note : \\
 \begin{center}
 $\begin{array}{cc} 
 \bar{R}_{01}= R_{01} (1:b,1:b)&\bar{\bar{R}}_{01}  = R_{01} (1:b,b+1:2b), \\
 \bar{R}_{11}  = R_{11} (1:b,1:b)&\bar{\bar{R}}_{11}  = R_{11} (1:b,b+1:2b). \\
 \end{array} 
 $
 \end{center}
 As in the previous level, we aim to eliminate b rows in each block, so we consider
 the matrix \\

 $$
 D_1 \Pi_1 =
 \left[
 \begin{array}{cc} 
 \left[ 
 \begin{array}{cccc} 
 I_b& & &  \\
  &I_{m/p - b}& & \\
 -D_{01}& &I_b& \\
  & & & I_{m/p - b} \\
 \end{array} 
 \right] & \\
  &\left[  \begin{array}{cccc} 
 I_b& & & \\
  &I_{m/p - b} & & \\
 -D_{11}&&I_b& \\
  & & &I_{m/p - b}  \\
 \end{array} 
 \right] \\
 \end{array} 
 \right] \times
 \left[ 
 \begin{array}{cc} 
 \Pi_{01}^T&\\
 &\Pi_{11}^T\\
 \end{array} 
 \right],
 $$

 where \\
 $$
 \begin{array}{c}  
 D_{01}=\bar{\bar{R}}_{01} ^T  (\bar{R}_{01} ^{-1})^T, \\
 D_{11}=\bar{\bar{R}}_{11} ^T  (\bar{R}_{01} ^{-1})^T. \\
 \end{array} 
 $$
 The matrices $ \bar{\Pi}_{01} $ and $ \bar{\Pi}_{11} $ are the permutations corresponding to
 the Strong RRQR factorizations of the two $2b \times b$ blocks. The matrices $\Pi_{01}$ and 
$\Pi_{11}$ can easily be deduced from the matrices $\bar{\Pi}_{01}$ and $\bar{\Pi}_{11}$ 
extended by the appropriate identity matrices to get matrices of size $2m/p\times 2m/p$.\\

 The multiplication of the block $D_0 \Pi_0 A(:,1:b)$ with the matrix
 $D_1 \Pi_1$ leeds to \\
 \begin{center}
 $$
 D_1 \Pi_1 D_0 \Pi_0 A(:,1:b) =
 \left[ 
 \begin{array}{c}
 \left[ \bar{\Pi}_{01}^T \times
 \left[ 
 \begin{array}{cc} 
 A_{01}   \\
 A_{11}\\
 \end{array} 
 \right] \right] (1:b,1:b)\\
 0_{2m/p - b}\\
 \\
 \left[ \bar{\Pi}_{11}^T \times
 \left[ 
 \begin{array}{cc} 
 A_{21}   \\
 A_{31}\\
 \end{array} 
 \right] \right] (1:b,1:b)\\
 0_{2m/p - b}\\
 \\
 \end{array} 
 \right].
 $$
 \end{center} 
 This matrix can be written as \\
 \begin{center} 
 $$
 D_1 \Pi_1 D_0 \Pi_0 A(:,1:b) =
 \left[ 
 \begin{array}{c}
 \left[ 
 \bar{\Pi}_{01}^T \times
 \left[ 
 \begin{array}{c} 
 (\Pi_{00}^T \times A_{0}) (1:b,1:b)\\
 (\Pi_{10}^T \times A_{1} )(1:b,1:b)\\
 \end{array} 
 \right] \right]  (1:b,1:b)\\
 0_{2m/p - b}\\
 \\
 \left[ 
 \bar{\Pi}_{11}^T \times 
 \left[ 
 \begin{array}{c}  
 (\Pi_{20}^T \times A_{2} (1:b,1:b))\\
 (\Pi_{30}^T \times A_{3} (1:b,1:b))\\
 \end{array} 
 \right] \right] (1:b,1:b)\\
 0_{2m/p - b}\\
 \\
 \end{array} 
 \right] =
 \left[
 \begin{array}{c}
 A_{02}  \\
 0_{2m/p - b}\\
 \\
 A_{12}\\
 0_{2m/p - b}\\
 \\
 \end{array} 
 \right]. 
 $$
 \end{center} 
 
 For the final step we consider the last $2b \times b$ block  \\

 $$
 \left[
 \begin{array}{c}
 A_{02}  \\
 A_{12}\\
 \end{array} 
 \right]. 
 $$

 We perform a Strong RRQR factorization on the transpose of this block and then we obtain  \\

 $$
 \left[ 
 \begin{array}{c}
 A_{02} \\
 A_{12}\
 \end{array} 
 \right] ^T \bar{\Pi}_{02} = Q_{02} R_{02}. \\
 $$
 We note  \\
 \begin{center}
 $\begin{array}{cc} 
 \bar{R}_{02} = R_{02} (1:b,1:b)&\bar{\bar{R}}_{02} = R_{02} (1:b,b+1:2b). \\
  \\
 \end{array} 
 $
 \end{center}

 We define the matrix $D_2 \Pi_2$ as \\
 $$
 D_2 \Pi_2 =
 \left[
 \begin{array}{cccccccc} 
 I_b& & & & & & & \\
  &I_{m/p - b}& & & & & & \\
 & &I_b& & & & & \\
 & & & I_{m/p - b}& & & &  \\
 -D_{02}& & & &I_b& & & \\
  & & & & & I_{m/p - b}& & \\
 & & & & & &I_b& \\
  & & & & & & &I_{m/p - b} \\
 \end{array} 
 \right] \times \Pi_{02}^T,
 $$
 where \\
 $$
 \begin{array}{c}  
 D_{02}=\bar{\bar{R}}_{02} ^T  (\bar{R}_{02}^{-1})^T,\\
 \\
 \end{array} 
 $$
 and the permutation matrix $\Pi_{02}$ can easily be deduced from the
 permutation matrix $\bar{\Pi}_{02}$. \\

 The multiplication of the block $D_1 \Pi_1 D_0 \Pi_0 A(:,1:b)$ by the matrix
 $D_2 \Pi_2$ leeds to : \\

 \begin{center}
 $$
 D_2 \Pi_2  D_1 \Pi_1 D_0 \Pi_0 A(:,1:b) =
 \left[ 
 \begin{array}{c}
  R_{021}^T  Q_{02}^T\\
  0_{4m/p - b} \\
 \end{array} 
 \right] =
 \Pi_{02}^T \times
 \left[ 
 \begin{array}{c}
 A_{02}\\
  0_{4m/p - b} \\
 \end{array} 
 \right]  =
 \left[ 
 \begin{array}{c}
 A_{03}\\
  0_{4m/p - b} \\
 \end{array} 
 \right]. 
 $$
 \end{center}

 If we consider the block $A(:,1:b)$ of the beginning and all the steps performed, we get :  \\
 $$
 D_2 \Pi_2 D_1 \Pi_1 D_0 \Pi_0 A(:,1:b) =
 \left[ 
 \begin{array}{c}
 \left[ 
 \bar{\Pi}_{02}^T \times
 \left[ 
 \begin{array}{c}
 A_{02}\\
 A_{12}\\
 \end{array} 
 \right] (1:b,1:b)\right] \\
 \\
 0_{4m/p - b}\\
 \end{array} 
 \right]
 $$
 We can also write \\
 $$
 D_2 \Pi_2 D_1 \Pi_1 D_0 \Pi_0 A(:,1:b) =
 \left[ 
 \begin{array}{c}
 \left[ 
 \bar{\Pi}_{02}^T \times
 \begin{array}{c}
 \left[ 
 \bar{\Pi}_{01}^T \times
 \left[ 
 \begin{array}{c} 
 (\Pi_{00}^T \times A_{0}) (1:b,1:b)\\
 \\
 (\Pi_{10}^T \times A_{1} )(1:b,1:b)\\
 \end{array} 
 \right] \right]  (1:b,1:b)\\
 \\
 \left[ 
 \bar{\Pi}_{11}^T \times
 \left[ 
 \begin{array}{c} 
 (\Pi_{20}^T \times A_{2}) (1:b,1:b)\\
 \\
 (\Pi_{30}^T \times A_{3} )(1:b,1:b)\\
 \end{array} 
 \right] \right]  (1:b,1:b)\\
 \end{array} 
 \right] (1:b,1:b)\\
 \\
 0_{4m/p - b}\\
 \end{array} 
 \right].
 $$

 We can also write  \\
 {\small
 $$
 A(:,1:b)=(D_2 \Pi_2 D_1 \Pi_1 D_0 \Pi_0)^{-1} 
 \left[ 
 \begin{array}{c}
 \left[ 
 \bar{\Pi}_{02}^T \times
 \begin{array}{c}
 \left[ 
 \bar{\Pi}_{01}^T \times
 \left[ 
 \begin{array}{c} 
 (\Pi_{00}^T \times A_{0}) (1:b,1:b)\\
 \\
 (\Pi_{10}^T \times A_{1} )(1:b,1:b)\\
 \end{array} 
 \right] \right]  (1:b,1:b)\\
 \\
 \left[ 
 \bar{\Pi}_{11}^T \times
 \left[ 
 \begin{array}{c} 
 (\Pi_{20}^T \times A_{2}) (1:b,1:b)\\
 \\
 (\Pi_{30}^T \times A_{3} )(1:b,1:b)\\
 \end{array} 
 \right] \right]  (1:b,1:b)\\
 \end{array} 
 \right] (1:b,1:b)\\
 \\
 0_{4m/p - b}\\
 \end{array} 
 \right].
 $$
 }

 Then, we have  \\
 {\small
 $$
 A(:,1:b)=\Pi_0^T D_0^{-1} \Pi_1^T D_1^{-1} \Pi_2^T D_2^{-1}
 \left[ 
 \begin{array}{c}
 \left[ 
 \bar{\Pi}_{02}^T \times
 \begin{array}{c}
 \left[ 
 \bar{\Pi}_{01}^T \times
 \left[ 
 \begin{array}{c} 
 (\Pi_{00}^T \times A_{0}) (1:b,1:b)\\
 \\
 (\Pi_{10}^T \times A_{1} )(1:b,1:b)\\
 \end{array} 
 \right] \right]  (1:b,1:b)\\
 \\
 \left[ 
 \bar{\Pi}_{11}^T \times
 \left[ 
 \begin{array}{c} 
 (\Pi_{20}^T \times A_{2}) (1:b,1:b)\\
 \\
 (\Pi_{30}^T \times A_{3} )(1:b,1:b)\\
 \end{array} 
 \right] \right]  (1:b,1:b)\\
 \end{array} 
 \right] (1:b,1:b)\\
 \\
 0_{4m/p - b}\\
 \end{array} 
 \right],
 $$
 }
 where \\
 $$
  \Pi_0^T =
  \left[
 \begin{array}{cccc} 
 \Pi_{00}& & &\\
 &\Pi_{10}& &\\
 & &\Pi_{20}&\\
 & & &\Pi_{30}\\
 \end{array} 
 \right],  
 $$
 $$
 D_0^{-1}=
 \left[
 \begin{array}{cccc} 
 \left[ 
 \begin{array}{cc} 
 I_b& \\
 D_{00}&I_{m/p - b}  \\
 \end{array} 
 \right] & & &\\
  &\left[  \begin{array}{cc} 
 I_b&  \\
 D_{10}&I_{m/p - b}  \\
 \end{array} 
 \right] & & \\
  & & \left[ \begin{array}{cc} 
 I_b& \\
 D_{20}&I_{m/p - b}  \\
 \end{array} 
 \right]  &  \\
 & & & \left[ 
 \begin{array}{cc} 
 I_b& \\
 D_{30}&I_{m/p - b}  \\
 \end{array} 
 \right] \\
 \end{array} 
 \right]. 
 $$
For each Strong RRQR performed previously, the corresponding trailing matrix is 
updated at the same time. Thus, at the end of the process, we have  \\
 $$
 A =  \Pi_0^T D_0^{-1} \Pi_1^T D_1^{-1} \Pi_2^T D_2^{-1}\times \left[
 \begin{array}{cc}
  A_{03}&\secA_{12}\\
 0_{4m/p - b}&\secA_{22}\\
 \end{array} 
 \right], 
 $$
where $ A_{03}$ is the $b\times b$ diagonal block containing the b selected pivot rows
 from the current panel. An additional GEPP should be performed on this block to get the
 full LU factorization.
$\secA_{22}$ is the trailing matrix already updated, and on which the block parallel
LU-PRRP should be continued.

\section*{Appendix F}
\label{sec:appendixF}
Here is the matlab code for generating the matrix $T$ used in Section 2
to define the generalized Wilkinson matrix on which GEPP fails. For any given 
integer $r > 0$, this generalized Wilkinson matrix is an upper triangular
semi-separable matrix with rank at most $r$ in all of its submatrices above  
the main diagonal. The entries are all negative and are chosen randomly. 
The Wilkinson matrix is for the special case where $r = 1$ and every entry 
above the main diagonal is 1.\\

\begin{verbatim}
	function [A] = counterexample_GEPP(n,r,u,v);
%       Function counterexample_GEPP generates a matrix which fails 
%	GEPP in terms of large element growth.
%       This is a generalization of the Wilkinson matrix.
%
	if (nargin == 2)
	      u = rand(n,r);
	      v = rand(n,r);
	      A = - triu(u * v');
	      for k = 2:n
	          umax = max(abs(A(k-1,k:n))) * (1 + 1/n);
	          A(k-1,k:n)  = A(k-1,k:n) / umax;
	      end
	      A = A - diag(diag(A)); 
	      A = A' + eye(n);
	      A(1:n-1,n) = ones(n-1,1);
	else
	      A = triu(u * v');
	      A = A - diag(diag(A));
	      A = A' + eye(n);
	      A(1:n-1,n) = ones(n-1,1);
	end

\end{verbatim}

\end{document}